\documentclass[oneside,12pt]{article}
\usepackage{graphicx}
\begin{document}
\title
{Handle additions producing essential surfaces}
\author{Ruifeng Qiu and Shicheng Wang
\thanks{Both authors are supported by  NSFC }}
\maketitle \pagestyle{myheadings} \markboth{R.F. Qiu and S.C.
Wang} {Handle additions producing essential surfaces}

\begin{abstract}
We construct a small, hyperbolic 3-manifold $M$ such that, for any
integer $g\geq 2$, there are infinitely many separating slopes $r$
in $\partial M$ so that $M(r)$, the 3-manifold obtained by
attaching a 2-handle to $M$ along $r$, is hyperbolic and contains
an essential separating closed surface of genus $g$.  The result
contrasts sharply with those known finiteness results on Dehn
filling, and it also contrasts sharply with the known finiteness
result on handle addition for the cases $g=0,1$. Our 3-manifold
$M$ is the complement of a hyperbolic, small knot in a handlebody
of genus 3.
\\keywords. Hyperbolic knot, Small knot, Handle addition.
\end{abstract}
\vskip 0.4 true cm

\begin{center}{\bf \S 1. Introduction.}\end{center}
\vskip 0.5 true cm

All terminologies not defined in the introduction are standard,
and will be defined later when they are needed.

All manifolds in this paper are orientable, all surfaces $F$ in
3-manifolds $M$ are embedded and proper, unless otherwise
specified. A surface $F\subset M$ is proper if $F\cap\partial
M=\partial F$.

Let $M$ be a compact 3-manifold. An incompressible,
$\partial$-incompressible surface $F$ in $M$ is  essential if it
is not parallel to $\partial M$. A 3-manifold $M$ is  simple if
$M$ is irreducible, $\partial$-irreducible, anannular and
atoroidal. In this paper, a compact 3-manifold $M$ is said to be
hyperbolic if $M$ with its toroidal boundary components removed
admits a complete hyperbolic structure with totally geodesic
boundary. By Thurston's theorem, a Haken 3-manifold is hyperbolic
if and only if it is simple. A knot $K$ in $M$ is hyperbolic if
$M_K$, the complement of $K$ in $M$, is hyperbolic. A 3-manifold
$M$ is {\it small} if $M$ contains no essential closed surface. A
knot $K$ in $M$ is {\it small} if $M_{K}$ is small.

A {\it slope} $r$ in $\partial M$ is an isotopy class of
unoriented essential simple closed curves in $F$.  We denote by
$M(r)$ the manifold obtained by attaching a 2-handle to $M$ along
a regular neighborhood of $r$ in $\partial M$ and then capping off
the possible spherical component with a 3-ball. In particular, if
$r$ lies in a toroidal component of $\partial M$, this operation
is known as Dehn filling.

Essential surface is a basic tool to study 3-manifolds and handle
addition is a basic method to construct 3-manifolds. A central
topic connecting  those two aspects is the following: \vskip 0.4
true cm

{\bf Question 1.} \ Let $M$ be a hyperbolic 3-manifold with
$\partial M\ne \emptyset$ which contains no essential closed
surface of genus $g$. How many slopes $r\subset
\partial M$ are there so that $M(r)$ contains an essential closed surface
of genus $g$?

{\bf Remark about Question 1.} \ The mapping class group of a
hyperbolic 3-manifold  is finite. The question is asked only for
hyperbolic 3-manifolds  to avoid possibly infinitely many slopes
produced by Dehn twist along essential discs or annuli.

The main result in this paper is the following:

{\bf Theorem 1.} \ There is a  small, hyperbolic  knot $K$ in a
handlebody $H$ of genus 3 such that, for any given integer $g\geq
2$, there are infinitely many separating slopes $r$ in $\partial
H$ such that $H_{K}(r)$ contains an essential separating closed
surface of genus $g$. Moreover those $H_K(r)$ are still
hyperbolic.

{\bf Comments.} \ Suppose $M$ is a hyperbolic 3-manifold with
non-empty $\partial M$.

(1) \ $\partial M$ is a torus. W. Thurston's pioneer result claims
that there are at most finitely many slopes on $\partial M$ such
that $M(r)$ are not hyperbolic, hence the number of slopes in
Question 1 is finite  when $g=0$, or $1$. Furthermore, the sharp
upper bounds of such slopes are given by Gordon and Lueck for
$g=0$, and by Gordon for $g=1$, see the survey paper [G].
Hatcher[Ha] proved that the number of slopes in Question 1 is
finite for any $g$.

(2) \ $\partial M$ has genus at least 2.  M. Scharlemann and Y-Q
Wu [SW] have shown that if $g=0$, or $1$ then there are only
finitely many separating slopes $r$ so that $M(r)$ contains an
essential closed surface of genus $g$. Very recently, M. Lackenby
[L] generalized Thurston's finiteness result to handlebody
attaching, that is to add 2-handles simultaneously. He proved
that, for a hyperbolic 3-manifold $M$, there is a finite set $C$
of exceptional curves on $\partial M$ so attaching a handlebody to
$M$ is still hyperbolike if none of those curves is attached to a
meridian disc of the handlebody.

(3) \ We [QW1] proved Theorem 1 for $g$ is even.

Theorem 1 and those finiteness results of [Th], [Ha], [SW] and [L]
give a global view about the answer of Question 1. In particular,
those finiteness results does not hold in general.

{\bf Outline of the proof of Theorem 1 and organization of the
paper.} In Section 2 we first construct a knot $K$ in the
handlebody $H$ of genus 3 for Theorem 1, then we construct
infinitely many surfaces $S_{g,l}$ of genus $g$ for each $g\ge 2$
such that (1) all those surfaces are disjoint from the given $K$,
therefore all $S_{g,l}\subset H_K$, (2)  for fixed $g$, all
$\partial S_{g,l}$ are connected and provide infinitely many
slopes in $\partial H$ when $l$ varies over from 1 to infinity.
Those $\partial S_{g,l}$ will be served as the slopes $r$ in
Theorem 1. We denote by $\hat S_{g,l}\subset H_K(\partial
S_{g,l})$ the closed surface of genus $g$ obtained by capping off
the boundary of $S_{g,l}$ with a disk. We will prove that $\hat
S_{g,l}$ is incompressible in $H_K(\partial S_{g,l})$ in Section
3. Sections 4 and 5 are devoted to prove that the knot $K$ is
hyperbolic and small.

\begin{center}{\bf \S 2. Construction of the knot $K$ and the surfaces
$S_{g,l}$ in $H$.}\end{center}

Let $H$ be a handlebody of genus 3. Suppose that $B_{1}, B_{2}$
and $B_{3}$ are basis disks of $H$, and $E_{1},E_{2}$ are two
separating disks in $H$ which separate $H$ into three solid tori
$J_{1}, J_{2}$ and $J_{3}$. See Figure 2.1.

Let $c$ be a closed curve in $\partial H$ as in Figure 2.2.  Then
$E_{1}\cup E_{2}$ separates $c$ into 10 arcs
$c_{1},\ldots,c_{10}$, where $c_{1},c_{3},c_{9}\subset J_{1}$ meet
$B_{1}$ in two, one, one points respectively; $c_{2}, c_{4}, c_6,
c_8, c_{10}\subset J_{2}$ meet $B_{2}$ in one, one, two, zero, one
points respectively; $c_{5}, c_{7}\subset J_{3}$ meet $B_{3}$ in
one, three points respectively.
\begin{center}
\includegraphics[totalheight=4cm]{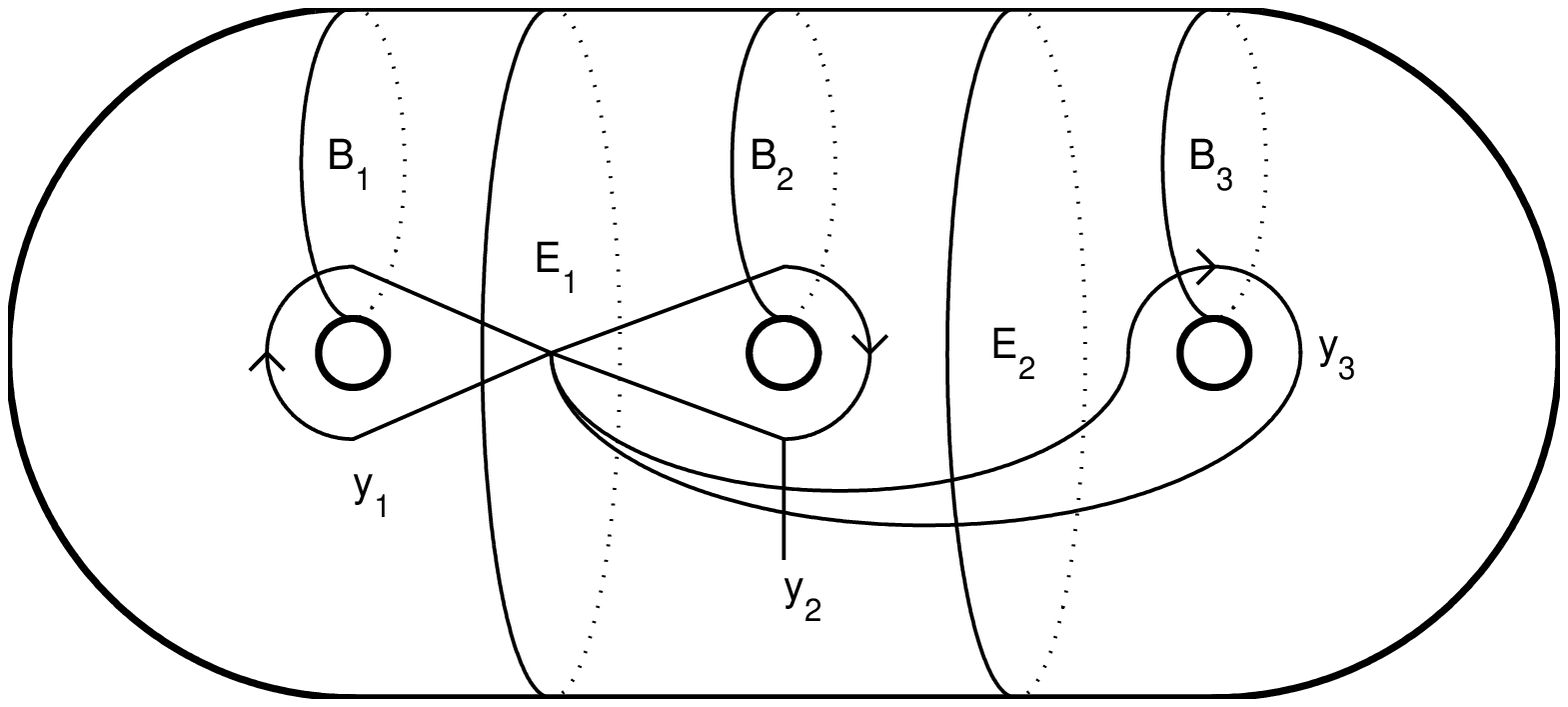}
\begin{center}
Figure 2.1

\end{center}
\end{center}

\begin{center}
\includegraphics[totalheight=5cm]{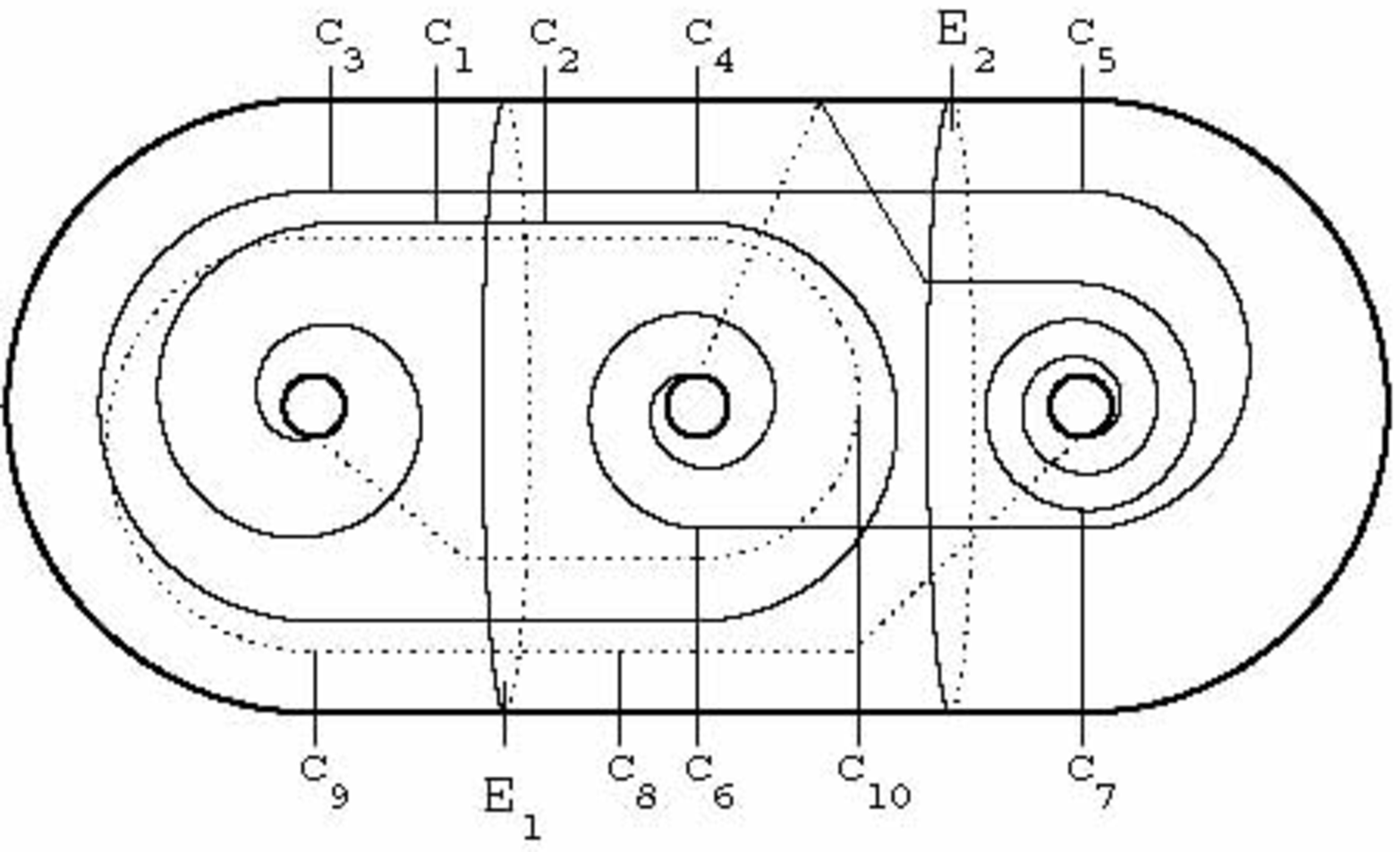}
\begin{center}
Figure 2.2
\end{center}
\end{center}

Let $u_{1},\ldots,u_{2g}, v_{1},\ldots,v_{2g}$ be $4g$ points
located on $\partial E_{1}$ in the cyclic order $u_1$, $u_3$, ...
, $u_{2i-1}$, ..., $u_{2g-1}$, $u_{2g}$, $u_{2g-2}$, ...
,$u_4$,$u_2$, $v_1$, $v_3$, ... , $v_{2i-1}$, ... $v_{2g-1}$,
$v_{2g}$, $v_{2g-2}$, ... ,$v_4$,$v_2$ as in Figure 2.3. By the
order of these points, $C$ can be isotoped so that $\partial c_1
=\{u_1, v_1\}$, $\partial c_2 =\{u_1, v_2\}$, $\partial c_{10}
=\{v_1, u_2\}$, $\partial c_3 =\{v_2, u_3\}$, $\partial c_{9}
=\{u_2, v_3\}$. Now suppose $\overline{u_{2i+1}v_{2i}}$ and
$\overline{v_{2i+1}u_{2i}}$ are arcs in $\partial J_{1}-intE_{1}$
parallel to $c_{3}$ and $c_{9}$ for $1\leq i\leq g-1$,
$\overline{u_{2}v_{1}}=c_{10}$, $\overline{v_{2}u_{1}}=c_{2}$,
$\overline{u_{2i}v_{2i-1}}$ and $\overline{v_{2i}u_{2i-1}}$ are
parallel arcs in $\partial (J_{2}\cup J_{3})-intE_{1}$, each of
which intersects $B_{2}$ in one point and $B_{3}$ in $l$ points as
in Figure 2.3 for $2\leq i\leq g$. Finally define
$\alpha_1$=$\overline{u_1v_1}$, and $\alpha_k$ is the union of
$\overline{v_{k-1}u_{k}}$, $\alpha_{k-1}$ and
$\overline{u_{k-1}v_{k}}$, $k=2,..., 2g$. Hence
$\alpha_{k-1}\subset \alpha_k$, $k=1,...., 2g$, is an increasing
sequence of arcs.

\begin{center}

\includegraphics[totalheight=5cm]{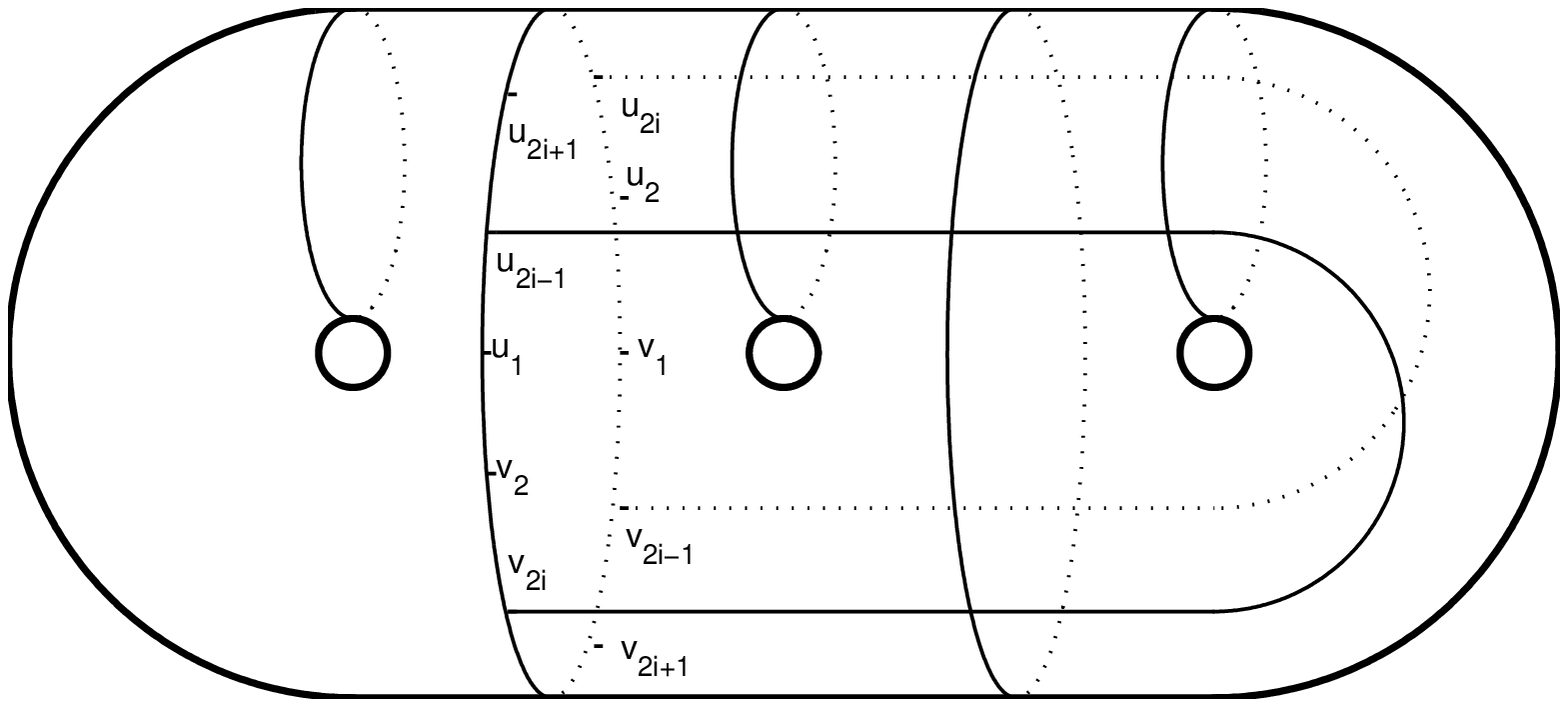}
\begin{center}
Figure 2.3
\end{center}
\end{center}

Let $\alpha\subset
\partial H$ be an arc which meets $\partial S$ exactly at its two
endpoints for a proper separating surfaces $S\subset H$. The
resulting surface of tubing $S$ along $\alpha$ in $H$, denoted by
$S(\alpha)$, is obtained by first attaching a 2-dimensional
1-handle $N(\alpha)\subset
\partial H$ to $S$, then making  the surface $S\cup N(\alpha)$ to
be proper, that is, pushing the  interior of $S\cup N(\alpha)$
into the interior of $H$. The image of $N(\alpha)$ after the
pushing is still denoted by $N(\alpha)$. In fact, $S\cup
N(\alpha)$ is a once punctured torus. Since $S$ is orientable and
separating, it is a direct observation that $S(\alpha)$ is still
separating and orientable.

Since $\alpha_1$ meets $E_1$ exactly in its two endpoints, we  do
tubing of $E_1$ along $\alpha_1$ to get a proper surface
$E_1(\alpha_1)$. Now $\alpha_2$ meets $E_1(\alpha_1)$ exactly in
its two endpoints, we do tubing of $E_1(\alpha_1)$ along
$\alpha_2$ to get $E_1(\alpha_1, \alpha_2)=E_1(\alpha_1)
(\alpha_2)$, where the tube $N(\alpha_2)$ is thinner and closer to
$\partial H$ so that it goes over the tube $N(\alpha_1)$. Hence
$E_1(\alpha_1, \alpha_2)$ is a properly embedded surface (indeed
an 1-punctured torus). By the same argument, we do tubing along
$\alpha_3,..., \alpha_{2g}$ to get a proper embedded surface
$E_1(\alpha_1,...., \alpha_{2g})$, denoted by $S_{g,l}$, in $H$,
moreover $S_{g,l}$ is orientable and separating.

Since the  surface $S_{g,l}$ is  obtained from the disc $E_1$ by
attaching $2g$ 1-handles to $E_1$ such that the ends of any two
handles are alternating,  $S_{g,l}$ is a once punctured orientable
surface of genus $g$. We summarize the facts discussed above as
the following:

{\bf Lemma 2.1.} \ $S_{g,l}$ is a once punctured surface of genus
$g$ and is separating in $H$.

Now let $K$ be a knot in $H$ obtained by pushing $C$ into $intH$
in the following way:

We first push $c_{6}$ into $intH$ deeply, and then we push
$C-c_{6}$ into $intH$ so that it stay between  $N(\alpha_3)$ and
$N(\alpha_4)$.

By observation, we have the following Lemma:

{\bf Lemma 2.2.} \ $K$ is disjoint from $S_{g,l}$ for all $g,l$.

 \vskip 0.5 true cm

\begin{center}{\bf \S 3. Proof of Theorem 1 by  assuming that $K$ is
hyperbolic and small.}\end{center}\vskip 0.5 true cm

We denote by $\hat S_{g,l}\subset H_K(\partial S_{g,l})\subset
H(\partial S_{g,l})$ the surface obtained by capping off the
boundary of $S_{g,l}$ with a disk. Then $\hat S_{g,l}$ is a closed
surface of genus $g$.

From the definition of $S_{g,l}$ for given genus $g$, $\partial
S_{g,l}$ provide infinitely many boundary slopes when  $l$ is
varied from 1 to infinity. Then Theorem 1 follows from the
following two propositions (the "moreover" part follows directly
from [SW]).

{\bf Proposition 3.0.} \ $K\subset H$ is a hyperbolic, small knot.

{\bf Proposition 3.1.} \  $\hat S_{g,l}$ is incompressible in
$H_K(\partial S_{g,l})$.

We will prove Proposition 3.1 in this section. Recall that a
surfaces $F$ in a 3-manifold is {\it compressible} if either $F$
is a 2-sphere which bounds a 3-ball, or there is an essential
simple closed curve in $F$ which bounds a disk in $M$; otherwise,
$F$ is {\it incompressible}. Hence Proposition 3.1 follows from
the following

{\bf Proposition 3.1*.} \  $\hat S_{g,l}$ is incompressible in
$H(\partial S_{g,l})$.

We choose the center of $E_1$ as the common base point for the
fundamental groups of $H$ and of all surfaces $S_{g,l}$.

Now $\pi_{1}(S_{g,l})$ is a free group of rank $2n$ generated by
$(x_{1},\ldots,x_{2n})$, where $x_i$ is the generator given by the
centerline of tube $N(\alpha_i)$; and $\pi_{1}(H)$ is a free group
of rank three generated by $y_{1}, y_{2}, y_{3}$ corresponding to
$B_{1}, B_{2}, B_{3}$ as in Figure 2.1.  Let $i:S_{g,l}\to H$ be
the inclusion. One can read $i_*(x_i)$ directly as the  words of
$y_1, y_2, y_3$ as follow:

$i_{\ast}(x_{1})=y_{1}^{2}$,

$i_{\ast}(x_{2})=y_{2}y_{1}^{2}y_{2}$,

$i_{\ast}(x_{3})=y_{1}y_{2}y_{1}^{2}y_{2}y_{1}$,

$i_{\ast}(x_{4})=y_{2}y_{3}^{l}y_{1}y_{2}y_{1}^{2}y_{2}y_{1}y_{2}y_{3}^{l}$,

and in general for $2 \le i\le g$,

$i_{\ast}(x_{2i-1})=y_{1}(y_{2}y_{3}^{l}y_{1})^{i-2}y_{2}y_{1}^{2}y_{2}(y_{1}y_{2}y_{3}^{l})^{i-2}y_{1}$.

$i_{\ast}(x_{2i})=(y_{2}y_{3}^{l}y_{1})^{i-1}y_{2}y_{1}^{2}y_{2}(y_{1}y_{2}y_{3}^{l})^{i-1}$.

{\bf Lemma 3.2.} \ $S_{g,l}$ is incompressible in $H$.

{\bf Proof.} \ The proof is the same as  that in [Q].\qquad Q.E.D.

Now $S_{g,l}$  separates $H$ into two components $P_{1}$ and
$P_{2}$ with $\partial P_{1}=T_{1}\cup S_{g,l}$ and $\partial
P_{2}=T_{2}\cup S_{g,l}$, where $T_1\cup T_2=\partial H$ and
$\partial T_1=\partial T_2=\partial S_{g,l}$.

{\bf Lemma 3.3.} \ $T_{i}$ is incompressible in $H$.

{\bf  Proof.} \ Note $H_{1}(H)=Z+Z+Z$ and with the three
generators  $y_{1}$, $y_{2}$ and $y_{3}$. By the above argument,
$i_{\ast} (H_{1}(S_{g,l}))$ is a subgroup of $H_{1}(H)$ generated
by $2y_{1}$, $2y_{2}$ and $2ly_{3}$. Thus $
H_{1}(H)/i_{\ast}(H_{1}(S_{g,l}))=Z_{2}\oplus Z_2\oplus Z_{2l}$ is
a finite group.

If $T_{i}$, $i=1$ or 2, is compressible. Then there is a
compressing disk $D_1$ in $H$ for $T_i$. Since $\partial
D\cap\partial S_{g,l}=\emptyset$ and $S_{g,l}$ is incompressible
in $H$, by a standard argument in 3-manifold topology, we may
assume that $D_1\cap S_{g,l}=\phi$. Furthermore, since $H$ is a
handlebody, we may also assume that $D_1$ is non-separating in
$H$. Thus there are two properly embedded disks $D_2$ and $D_3$ in
$H$ such that $\bigl\{D_1,D_2, D_3\bigr\}$ is a set of basis disks
of $H$. Let $z_{1}$, $z_{2}$ and $z_{3}$ be generators of
$\pi_{1}(H)$ corresponding to $D_1, D_2$ and $D_3$. Since
$S_{g,l}$ misses $D_1$, $i_{\ast}(\pi_{1}(S_{g,l}))\subset G$
where $G$ is a subgroup of  $\pi_{1}(H)$ generated by $z_{2}$ and
$z_{3}$. Then clearly $H_{1}(H)/i_*(H_{1}(S_{g,l}))$ is an
infinite group, a contradiction.\qquad Q.E.D.

{\bf Proof of Proposition 3.1*.} \ Since $H$ is a handlebody and
$S_{g,l}$ is incompressible in $H$,  $P_{1}$ and $P_{2}$ are
handlebodies. By Lemmas 3.2, 3.3 and the Handle Addition
lemma[Jo], $\hat S_{g,l}$ is incompressible in $P_{i}(\partial
S_{g,l})$ for $i=1,2.$ Since $H(\partial S_{g,l})=P_{1}(\partial
S_{g,l})\cup_{\hat S_{g,l}} P_{2}(\partial S_{g,l})$, $\hat
S_{g,l}$ is incompressible in $H(\partial S_{g,l})$.\qquad Q.E.D.

\begin{center}{\bf \S 4. $H_k$ is irreducible, $\partial$-irreducible,
anannular}\end{center} \vskip 0.5 true cm

By the construction, $K$ is cut by $E_{1}\cup E_{2}$ into ten arcs
$a_{1},\ldots, a_{10}$ where $a_{i}$ is obtained by pushing
$c_{i}$ into $intH$. Now let $N(K)=K\times D$ be a regular
neighborhood of $K$ in $H$ such that the product structure has
been adjusted so that $\cup_{i=1}^{10}\partial a_i\times D\subset
E_1\cup E_2$. Let $H_{K}=H-intN(K)$ and $F_{i}=E_{i}-intN(K)$. We
denote by $M_{i}=H_{K}\cap J_{i}$, $i=1,2,3$, and $T=\partial
(K\times D)$. Then $F_{1}\cup F_{2}$ separates $T$ into ten annuli
$A_{1},\ldots, A_{10}$ such that $A_{i}=a_{i}\times \partial D$.

Moreover $K$ and $C$ bound a non-embedded annulus $A_*$ which is
cut by $E_{1}\cup E_{2}$ into ten disk $D_{1*},..., D_{10*}$ in
$H$. Note that $D_*=\cup_{i\ne 6}D_{i*}$ is still a  disk. Let
$D_i=D_{i*}\cap H_K$ for $i\neq 6$. Then $D_i$ is a proper disk in
some $M_l$ and $\cup_{i\ne 6}D_{i}$ is a still disk, see Figure
4.1. Now we number $\partial A_{i}$ such that $\partial_{1}
A_{i}=\partial_{2} A_{i-1}$ and $\partial_{2} A_{i}=\partial_{1}
A_{i+1}$. For $i\neq 6$, let $W_{i}=\overline{\partial N(D_i\cup
A_i)-\partial M_l} $. Then $W_{i}$ is a proper separating disk in
$M_{l}$. Each $W_{i}$ intersects $F_{1}\cup F_{2}$ in two arcs
$l_{i}$ and $l_{i+1}$. Note that $W=\cup_{i\neq 6}W_{i}$ is a
disk. Thus $\partial W$ is a union of two arcs in $\partial H$ and
$ l_{6}\cup l_{7}$, see Figure 4.1. Since $c_{3}, c_{9}$ are
parallel in $\partial J_{1}-intE_{1}$, there are two arcs parallel
to $c_{3}$ in $\partial J_{1}-intE_{1}$, say $l^{'}, l^{''}$, and
two arcs in $F_{1}$, say $l^{1}, l^{2}$, such that $l^{'}\cup
l^{''}\cup l^{1}\cup l^{2}$ bounds a disk $W^{'}$ which separates
$M_{1}$ into two handlebodies of genus two $H^{1}, H^{2}$ with
$A_{1}\subset H^{1}$ and $A_{3}, A_{9}\subset H^{2}$. We denote by
$\mu$ the meridian slope on $T$ and $\tau$ the longitude slope on
$T$.
\begin{center}
\includegraphics[totalheight=6cm]{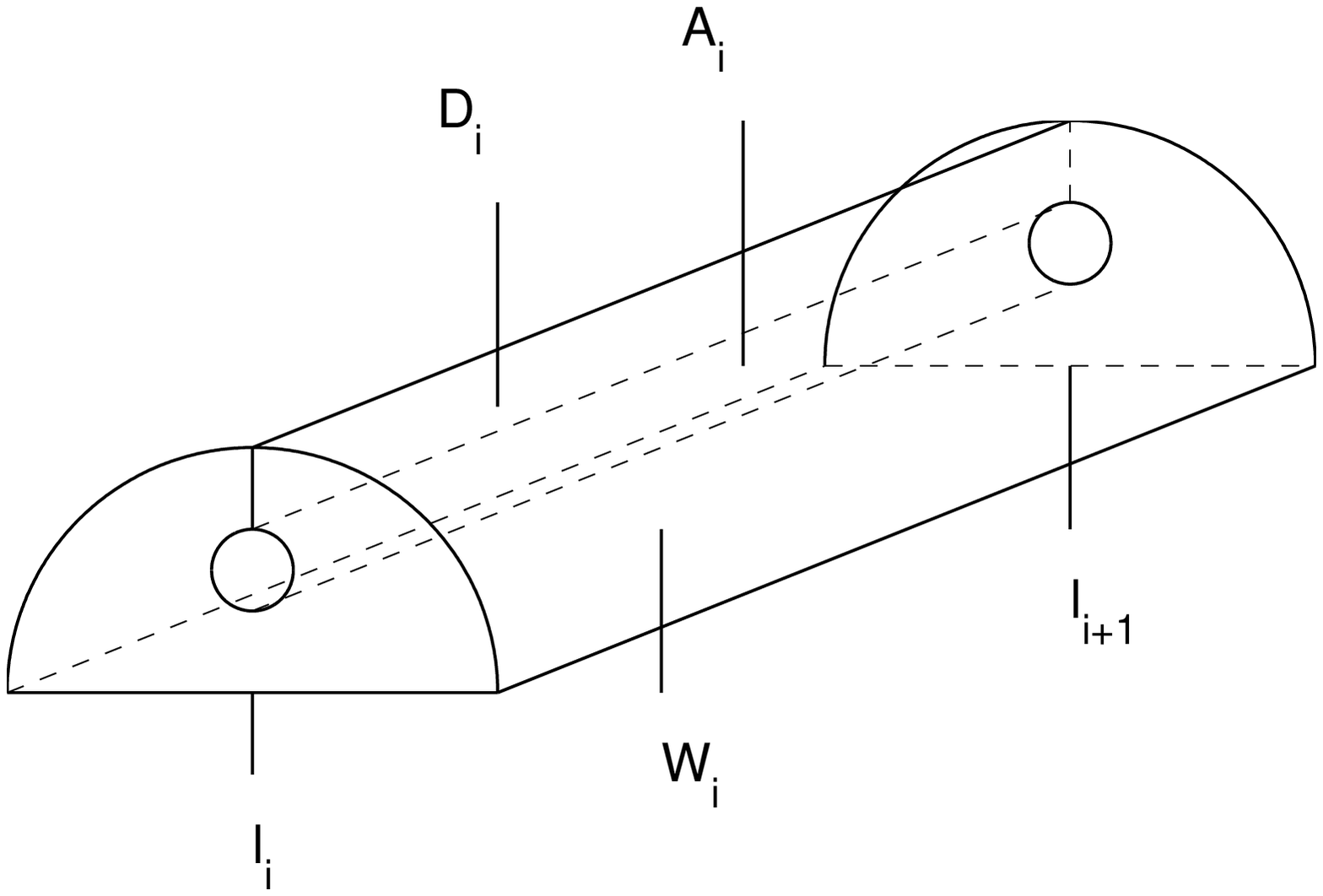}
\begin{center}
Figure 4.1
\end{center}
\end{center}

We list some elementary facts about $K$ and $a_i$ as the
following:

{\bf Lemma 4.0.}

(1) $K\ne 1$ in $\pi_1(H)$.

(2) Suppose $a_i\subset J_m$, $i\ne 4, 8$. Let $b_i\subset E_1\cup
E_2$ be a given arc with $\partial b_i=\partial a_i$ and $B\subset
J_m$ be a non-separating proper disk. Then $a_{i}\cup b_{i}$
intersects $\partial B$ in at least one point for all $i$, in at
least three points when $i= 7$,  in at least two points when $i=1,
6.$

(3) There is no relative homotopy on $(J_m, E_1\cup E_2)$ which
sends $a_i$ to $E_1\cup E_2$.

Recall that a 3-manifold $M$ is {\it irreducible} if it contains
no essential 2-spheres, $M$ is {\it $\partial$-irreducible} if
$\partial M$ is incompressible, $M$ is atoroidal if it contains no
essential tori, $M$ is {\it anannular} if it contains no essential
annuli.

{\bf Lemma  4.1.} \ $H_{K}$ is irreducible.

{\bf Proof.} \ Suppose that $H_{K}$ is reducible.  Then there is
an essential 2-sphere $S$ in $H_{K}$.  Since $H$ is irreducible,
$S$ bounds a 3-ball $B^3$ in $H$ and  $K\subset B^3$, which
contradicts (1) of Lemma 4.0.\qquad Q.E.D.

Recall that $F$ is {\it $\partial$-compressible} if there is an
essential arc $a$ in $F$ which, together with an  arc $b$ in
$\partial M$, bounds a disk $D$ in $M$ such that $D\cap F=a$;
otherwise, $F$ is $\partial$-incompressible.

{\bf Lemma 4.2.} \ $F_{1}\cup F_2$ is incompressible and
$\partial$-incompressible in $H_{K}$.

{\bf Proof.} \  Suppose first that $F_{1}\cup F_{2}$  is
compressible in $H_{K}$. Then there is a disk $B$ in $M$ such that
$B\cap (F_{1}\cup F_{2})=\partial B$ and $\partial B$ is an
essential circle on $F_{1}\cup F_{2}$. Without loss of generality,
we assume that $\partial B\subset F_{1}$ and $B\subset M_{2}$.
Denote by $B^{'}$ the disk bounded by $\partial B$ in $E_{1}$.
Then $B\cup B^{'}$ is a 2-sphere $S\subset J_2$. Then it follows
easily from Lemma 4.1 that $S$ bounds a 3-ball $B^3$ in $J_{2}$.
Since $\partial B$ is essential in $F_{1}$, $B^{'}$ contains at
least one component of $\partial a_i$. Since $S$ is separating and
$a_i$ is connected, we must have $(a_{i}, \partial a_i)\subset
(B^3, B')$, which provides a relative homotopy on $(J_2, E_1)$
which sends $a_i$ to $E_1$. It contradicts (2) of Lemma 4.0.

Suppose now that $F_{1}\cup F_{2}$ is $\partial$-compressible in
$H$. Then there is an essential arc $a$ in $F_{1}\cup F_{2}$
which, with an arc $b$ in $\partial H_{k}$, bounds a disk $B$ in
$H_{K}$ with $B\cap (F_{1}\cup F_{2})=a$. Without loss of
generality, we assume that $a\subset F_{2}$ and $B\subset M_{2}$.
There are two cases:

1) \ $b\subset T$. Then $b$ is a proper arc in one of  $A_{4}$,
$A_{6}$ $A_{8}$, say $A_{6}$. If $b$ is not essential in $A_6$,
then $a$ and an arc $b'$ in $\partial A_6$ forms an essential
circle in $F_2$ which bounds a disc in $M_2$. This contradicts the
incompressibility of $F_2$ we just proved. If $b$  is essential in
$A_6$, then the disk $B$ provides a relative homotopy on $(J_2,
E_2)$ sending $a_6$ to $E_2$, which contradicts (2) of Lemma 4.0.

2) \ $b\subset \partial H$. If $B$ is non-separating in $J_{2}$,
then $b_{6}$ can be chosen so that $a_{6}\cup b_{6}$ intersects
$\partial B$ in at most one point where $b_{6}$ is an arc in
$E_{2}$ connecting the two endpoints of $a_{6}$, contradicting (2)
of Lemma 4.0.  If $B$ is separating in $J_{2}$, then $B$ separates
$J_{2}$ into a 3-ball $B^3$ and a solid torus $J$. We denote by
$D_{1}, D_{2}$ the two components of $E_{2}-a$. Since $a$ is
essential in $F_{2}$, each of $intD_{1}$ and $intD_{2}$ contains
at least one endpoints of $a_{4}$, $a_{6}$ and $a_{8}$.

Suppose that $D_{1}\subset B^3$ and $D_{2}\cup E_{1}\subset J$. By
the construction, $\partial_{1} a_{4}, \partial_{1} a_{8}\subset
E_{1}$, $\partial_{2} a_{4}, \partial_{2} a_{8}\subset E_{2}$, and
$\partial a_{6}\subset E_{2}$. Since $a_{4}, a_{6}$ and $a_{8}$
are disjoint from $B$, $a_{4}, a_{8}\subset J$ and $a_{6}\subset
B^{3}$. It contradicts (2) of Lemma 4.0.

Suppose that $D_{1}\subset J$ and $D_{2}\cup E_{1}\subset B^{3}$.
Then $a_{2},a_{10}\subset B^{3}$. It contradicts (2) of Lemma
4.0.\qquad Q.E.D.

{\bf Lemma 4.3.} \ $H_{K}$ is $\partial$-irreducible.

{\bf Proof.} \ Suppose $H_{K}$ is $\partial$-reducible. Let $B$ be
a compressing disk of $\partial H_{K}$. If $\partial B\subset T$,
then $H_K$ contains an essential 2-sphere, which contradicts Lemma
4.1. Below we assume that $\partial B\subset
\partial H$. Since $F_{1}\cup F_{2}$ is incompressible and
$\partial$-incompressible in $H_{K}$ (Lemma 4.2), by a standard
cut and paste argument, we may assume that $B\cap (F_{1} \cup
F_{2})=\emptyset$. We assume that $B\subset M_{2}$. (The other
cases are similar.) Then  $B$ misses $b_6$.  If $B$ is
non-separating in $J_{2}$,   By (2) of Lemma 4.0, $B$ intersects
$a_6$, a contradiction. If $B$ is separating, then $B$ separates a
3-ball $B^3$ from $J_2$. Since $\partial B$ is essential in
$\partial H_K$, there are two cases: Either $ B^3$ contains only
one of $ E_{1}$ and $ E_{2}$, say $E_1$, then $a_{8}\cap B\neq
\emptyset$, a contradiction; or $B^3$ contains both $E_{1}$ and
$E_{2}$, then there is a relative homotopy on $(J_2, E_2)$ sending
$a_6$ to $E_2$, which contradicts (2) of Lemma 4.0.\qquad Q.E.D.

{\bf Lemma 4.4.} \ $M$ is anannular.

{\bf Proof.} \ Suppose $H_{K}$ contains an essential annulus $A$.
Assume that

(**) $| A\cap (F_{1}\cup F_{2})|$ is minimal among all essential
annuli in $H_{K}$.

By Lemma 4.2, (**) and the proof of Lemma 4.3, each component of
$A\cap (F_{1}\cup F_{2})$ is  essential in both $A$ and
$(F_{1}\cup F_{2})$. There are three cases:

1. \ $\partial A\subset T$.  Now $A$ is separating in $H_{k}$;
otherwise, $H$ contains either a non-separating 2-sphere or a
non-separating tori. Hence the union of $A$ and an annulus $A^{'}$
on $T$ makes a separating torus $T^{'}$, cutting off a manifold
with boundary $T\cup T^{'}$. Since $M$ is irreducible, $T^{'}$ is
incompressible, so by Lemma 5.5 $T^{'}$ is parallel to $T$, which
implies that $A$ is inessential. (The arguments in Section 5 are
independent of the arguments in Section 4.)

2. $\partial_{1} A\subset T$ and $\partial_{2}A\subset \partial
H$.

By Lemma 4.3, both $\partial H$ and $T$ are incompressible  in
$H_K$. Clearly $H_K$ is not homeomorphic to $T\times I$. Since
both Dehn fillings along $\mu$ and $\partial A_1$ compress
$\partial H$, by an important theorem in Dehn filling, $\Delta
(\partial_{1} A,\mu)\leq 1$. [See 2.4.3 CGLS.]

We first suppose that $\partial_{1} A$ is the meridian slope
$\mu$. Then $\partial_{1} A$ is disjoint from $F_{1}\cup F_{2}$.
Now we claim that $A$ is disjoint from $F_{1}\cup F_{2}$.

Suppose, otherwise,  that $A\cap (F_{1}\cup F_{2})\neq\emptyset$.
Since $F_{1}\cup F_{2}$ is incompressible and
$\partial$-incompressible in $H_{K}$ (Lemma 4.2), by a standard
cut and paste argument, we may assume that $\partial_{2} A\cap
(F_{1} \cup F_{2})=\emptyset$. Now each component of $A\cap
(F_{1}\cup F_{2})$ is an essential simple closed curve in $A$. Let
$a$ be an outermost circle in $A\cap (F_{1}\cup F_{2})$. Then $a$
and $\partial_1 A$ bound an annulus $A^*$ in $A$ such that
$intA^*$ is disjoint from $F_{1}\cup F_{2}$. We may assume that
$a\subset F_{1}$, and $\partial_{1} A \subset A_{i}$ for some $i$.
Let $B^{*}$ be the disk bounded by $a$ on $E_{1}$ and $D$ be the
meridian disk of $N(K)$ bounded by $\partial_1 A$. Since $a$ is
essential on $F_{1}$, $B^{*}$ contains at least one component of
$\partial F_{1}$. In $H$, $B^{*}\cup A^*\cup D$ is a separating
2-sphere $S^{2}$ which bounds a 3-ball $B^3$. For $j\neq i$, if
$\partial_{1} a_{j}\subset B^{*}$, then $\partial_{2} a_{j}\subset
B^{*}$ and $a_j\subset B^3$. This possibility is ruled by (2) of
Lemma 4.0. Note also that $\partial_{1} a_{i}\subset B^{*}$ and
$\partial_{2} a_{i}$ is not contained in $B^{*}$. Now we denote by
$A^{'}$ the annulus bounded by $a$ and $\partial_{1}
a_{i}\times\partial D=\partial_{1} A_{i}$ in $F_{1}$. Then
$A^{*}\cup A^{'}$ is isotopic to an annulus disjoint from
$F_{1}\cup F_{2}$. By the above argument, $A^{*}\cup A^{'}$ is
inessential. Thus we can properly isotope $A$ by pushing the
annulus $A^{*}$ to the other side of $F_1$ to reduce $|A\cap
(F_{1}\cup F_{2})|$, which contradicts (**).

We may assume that $A$ is contained in $M_{2}$. Let $D$ be the
meridian disk of $N(K)$ bounded by $\partial_1 A$ and
$B=A\cup_{\partial_{1} A} D$. Then $B$ is a proper disk in $J_2$,
meeting $K$ in exactly one point, hence $B$ is a meridian disk of
$J_{2}$. Let $b_{6}$ be an arc on $E_{2}$ connecting the two
endpoints of $c_{6}$. Then $c_{6}\cup b_{6}$ would be a closed
curve of winding number 2 in the solid torus $J_{2}$ intersecting
$B$ at most once, which is absurd.







Now we suppose that $\Delta (\partial_{1} A, \mu)=1$.

Then  $A$ is cut by $ (F_{1}\cup F_{2})$ into ten squares $S_i$,
$i=1,...,10$, each $S_i$ has two opposite sides in $F_1\cup F_2$
and remaining two sides the longitude arc $a_i$ in $A_{i}$ and
$a_{i}^{*}\subset
\partial H$. Let $b^*_{2}$ be the arc connecting the two endpoints
of $a_{2}^{*}$ in $E_{1}$ and $b^*_{6}$ be the arc connecting the
two endpoints of $a_{6}^{*}$ in $E_{2}$. Then the two simple
closed curves $b^*_{2}\cup a_{2}^{*}$ and $b^*_{6}\cup a_{6}^{*}$
on $\partial J_2$ are disjoint. But in $\pi_{1}(J_{2})$,
$b^*_{2}\cup a_{2}^{*}=y_{2}$ and $b^*_{6}\cup
a_{6}^{*}=y_{2}^{2}$, a contradiction.

3. $\partial A\subset \partial H$.

Suppose first that $A\cap (F_{1}\cup F_{2})=\emptyset$. Then $A$
is contained in one of $M_{1}$, $M_{2}$ and $M_{3}$. Since $A$ is
essential and $H_{K}$ is $\partial$-irreducible, $A$ is disjoint
from $D_{i}$ for $i\neq 6$. Since each component of $\partial
H\cap  J_{1}-c_{1}\cup c_{3}$ and $\partial H\cap J_{3}-c_{5}\cup
c_{7}$ is a disc,  $A\subset M_{2}$. Since $A$ is disjoint from
$c_{2}, c_{4}, c_{8}, c_{10}$, each component of $\partial A$
intersects $B_{2}$ in only one point in $J_{2}$(see Figure 2.2).
Thus $A$ is isotopic to each component of $\partial J_{2}-\partial
A$ in $J_{2}$. It means that $A$ is not essential in $M_{2}$, a
contradiction.

Now suppose that $A\cap (F_{1}\cup F_{2})\neq \emptyset$.  There
are two subcases:

1) Each component of $A\cap(F_{1}\cup F_{2})$ is an essential
circle. Now let $a$ be an outermost component of $A\cap(F_{1}\cup
F_{2})$. That means that $\partial_{1} A$, with $a$, bounds an
annulus $A^{*}$ in $A$ such that $A^{*}\cap (F_{1}\cup F_{2})=a$.
Then $A^{*}\subset M_{i}$. We denote by $B^{*}$ the disk bounded
by $a$ in $E_{1}\cup E_{2}$. Let $D^{*}=A^{*}\cup B^{*}$. Then
$D^{*}$ be a disk. Let $D$ be the disk obtained from $D^{*}$ by
pushing $B^{*}$ slightly into $J_{l}$. Then $D$ is a properly
embedding disk in $J_{l}$ such that $D$ intersects each $a_{i}$ in
at most two points. Furthermore, if $D$ intersects $a_{i}$ in two
points for some $i$, then the two endpoints of $a_{i}$ lie in
$B^{*}$. Thus, in this case,  the algebraic intersection number of
$a_{i}$ and $D$ is $0$. Now by Lemma 4.0, $A^{*}$ is separating in
$J_{l}$.

Suppose that $A^{*}$ is contained in one of $J_{1}$ and $J_{3}$,
say $J_{1}$. Then $\partial_{1} A$ is parallel to $\partial
E_{1}$. We denote by $A^{'}$ the annulus bounded by $\partial_{1}
A$ and $a$ in $\partial J_{1}$. Since $a$ is essential in $F_{1}$,
$B_{*}$ contains at least one endpoint of $a_{1},a_{3},a_{9}$.
Furthermore, $\partial_{1} a_{i}\subset B^{*}$ if and only if
$\partial_{2} a_{i}\subset B^{*}$. Now if, for some $j$,
$\partial_{1} a_{j}\subset A^{'}$, then $\partial_{2} a_{j}\subset
A^{'}$. That means that $a_{j}$ is disjoint from $B_{1}$ as in
Figure 2.1, a contradiction. Thus for each $i,j$, $\partial_{j}
a_{i}\subset B^{*}$. That means that $a$ is parallel to $\partial
E_{1}$ in $F_{1}$. Now $\partial D_{i}$ intersects  each component
of $\partial A^{*}$ in two points for $i=1,3,9.$ That means that
$D_{i}$ intersects $A^{*}$ in two arcs such that each of the two
arcs has its two endpoints in distinct components of $\partial
A^{*}$. (Otherwise, since $\partial_{1} A$ is isotopic to
$\partial E_{1}$, $a_{i}\cup b_{i}=1$ in $\pi_{1}(J_{1})$ where
$b_{i}$ is an arc in $\partial E_{1}$ connecting the two endpoints
of $a_{i}$, a contradiction.) Thus we can push $\partial_{1} A$
into $M_{2}$ to reduce $|A\cap (F_{1}\cup F_{2})|$.

Suppose now that $A^{*}\subset M_{2}$. Without loss of generality,
we assume that $a\subset F_{1}$. We denote by $A^{'}$ the annulus
bounded by $\partial E_{1}$ and $a$ in $E_{1}$. Then $A^{'}$ and
$B^{*}$ lie in distinct sides of $J_{2}-A^{*}$. If $\partial_{1}
A$ is isotopic to $\partial E_{2}$, then $a_{6}\cup b_{6}=1$ in
$\pi_{1}(J_{2})$ where $b_{6}$ is an arc in $E_{2}$ connecting the
two endpoints of $a_{6}$ , a contradiction. If $\partial_{1} A$
bounds a disk $D$ in $\partial J_{2}$ such that $E_{1},
E_{2}\subset D$, then $a_{4}\cup a_{8}\cup b^{1}\cup b^{2}=1$ in
$\pi_{1}(J_{2})$ where $b^{i}$ is an arc in $E_{i}$ connecting the
endpoints of $a_{4i}$ and $a_{8}$, a contradiction. Now
$\partial_{1} A$ is isotopic to $\partial E_{1}$. Then $D_{4}$
intersects $A^{*}$ in an arc. By the above argument, we can push
$\partial_{1} A$ into $M_{1}$ to reduce $|A\cap (F_{1}\cup
F_{2})|$.
\begin{center}
\includegraphics[totalheight=5cm]{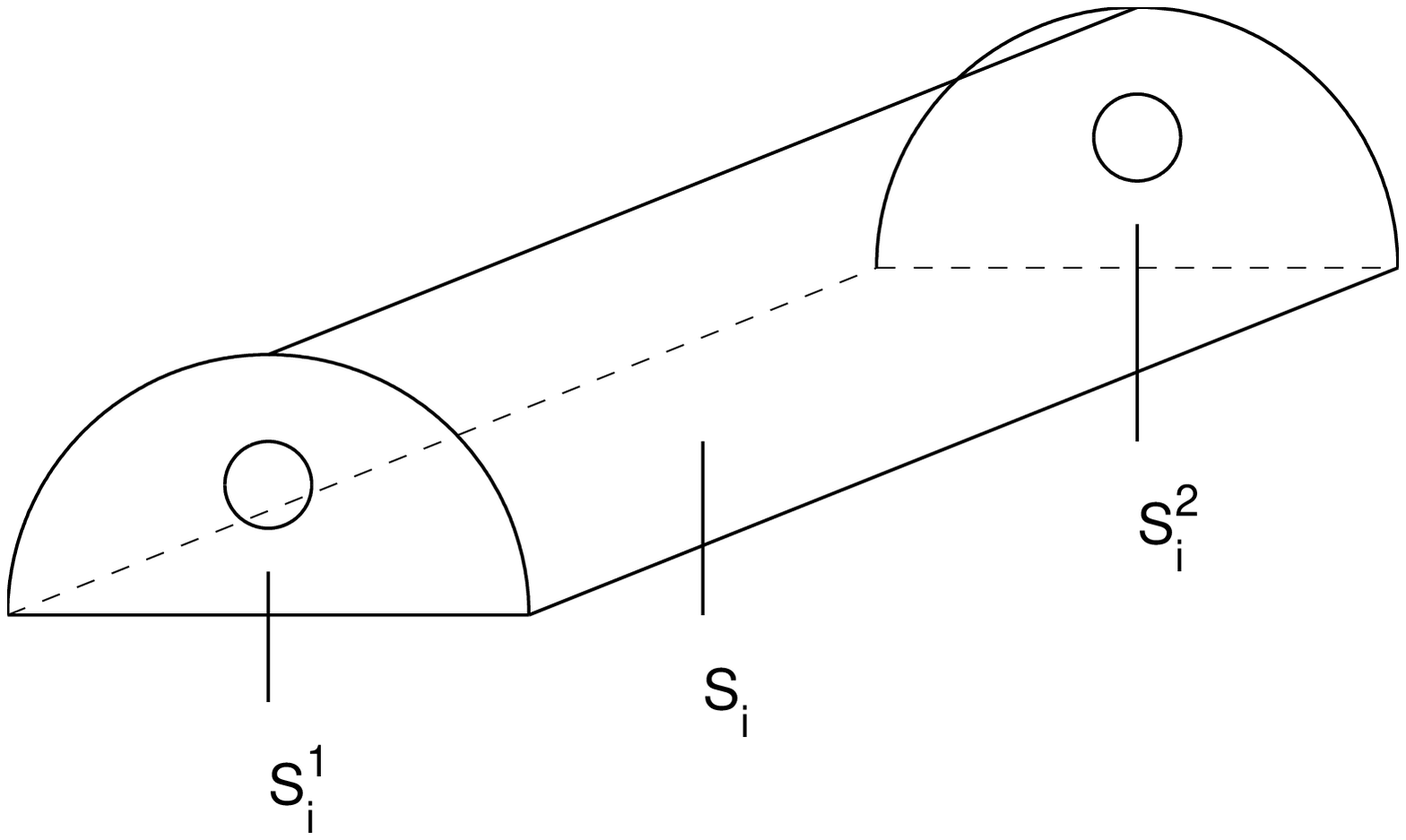}
\begin{center}
Figure 4.2
\end{center}
\end{center}

2) Each component of $A\cap(F_{1}\cup F_{2})$ is an essential arc.
Then $F_{1}\cup F_{2}$ cuts $A$ into proper squares $S_{i}$ in
$M_l\subset J_l$, each $S_i$ having two opposite sides in $F_1\cup
F_2$ and the remaining two sides  in $\partial H$. If
$S_{i}\subset J_{l}$ for $l=2$ or $3$, then $S_i$ is a separating
disc in $J_l$. Otherwise, say $S_i$ is a non-separating disc in
$J_2$, by the same reason as that at the end of the proof of Lemma
4.3, $S\cap (F_1\cup F_2)$ consisting of  two proper arcs in
$E_1\cup E_2$ implies that $b_6$ can be chosen so that $b_{6}$
intersects $\partial S_i$ in at most two points; furthermore, if
$b_{6}$ intersects $\partial S_i$ in two points then $S_{i}\cap
F_{1}=\emptyset$ and $S_{i}\cap b_{2}=\emptyset$, where $b_{i}$ is
an arc in $E_{1}\cup E_{2}$ connecting the two endpoints of
$a_{i}$. This means that $S_i$ meets $a_{2}$ or $a_6$ by Lemma 4.0
(1), a contradiction. Now each $ S_{i}$ cuts off a 3-ball
$B^{3}_{i}$ from $J_l$  for $l=2$ or $3$ as in Figure 4.2. Let
$S_{i}^{1}$ and $S_{i}^{2}$ be the two disks of $B^3_i\cap
(E_1\cup E_2)$ and $S_{i}\subset J_{l}$ where $l=2$ or $3$. By (2)
of Lemma 4.0, we have

i) \ $\partial_1 a_{j} \subset S_{i}^{1}$ if and only if
$\partial_2 a_{j}\subset S_{i}^{2}$.

ii) If  $a_{j}$ is contained in $B^{3}_{i}$, then $a_{l}$ is not
contained in $B^{3}_{i}$.

This means that for each $i$, there is only one boundary component
of $F_{1}\cup F_{2}$ lying in each of $S_{i}^{1}$ and $S_{i}^{2}$.
Thus if for some $i$, $S_{i}$ lies in $M_{1}$, then $S_{i}$ is
also separating in $J_{1}$. Otherwise, say $S_{i}$ is
non-separating in $J_{1}$. Then, by i) and ii), the three circles
$a_{1}\cup b_{1}, a_{3}\cup b_{3}, a_{9}\cup b_{9}$ intersect
$S_{i}$ in two points, a contradiction. It follows that $S_{i}$ is
also as in Figure 4.2 and $A$ cuts off a solid torus $P$ from $H$.
Thus $D_{i*}$ can be chosen to be disjoint from $A$ even if $i=6$.
This means that $K$  and  a component of $\partial A$ bounds an
annulus, which has been ruled out in Case 2.\qquad Q.E.D.

\vskip 0.5 true cm

\begin{center}{\bf  \S 5. $H_{K}$ contains no closed essential
surface}\end{center}\vskip 0.5 true cm

Suppose $H_{K}$ contains  essential closed surfaces. Let $W,
W^{'}$ and $W_{i}$ be the disk defined in Section 4. Denote by
$X(F)$ the union of the components of $F\cap M_{1}$ isotopic to
$\partial H\cap M_{1} $. We define the complexity on the essential
closed surfaces $F$ in $H_K$ by the following quadruple, in
lexicographic order.

$$C(F)=(|F\cap W|, |F\cap F_{2}|, |(F\cap M_{1}-X(F))
\cap W^{'}|,|F\cap F_{1}|).$$ Suppose $F$ realizes the minimality
of $C(F)$. By the minimality of $C(F)$ and the standard argument
in 3-manifold topology, we have Lemma 5.0.

\begin{center}
\includegraphics[totalheight=4cm]{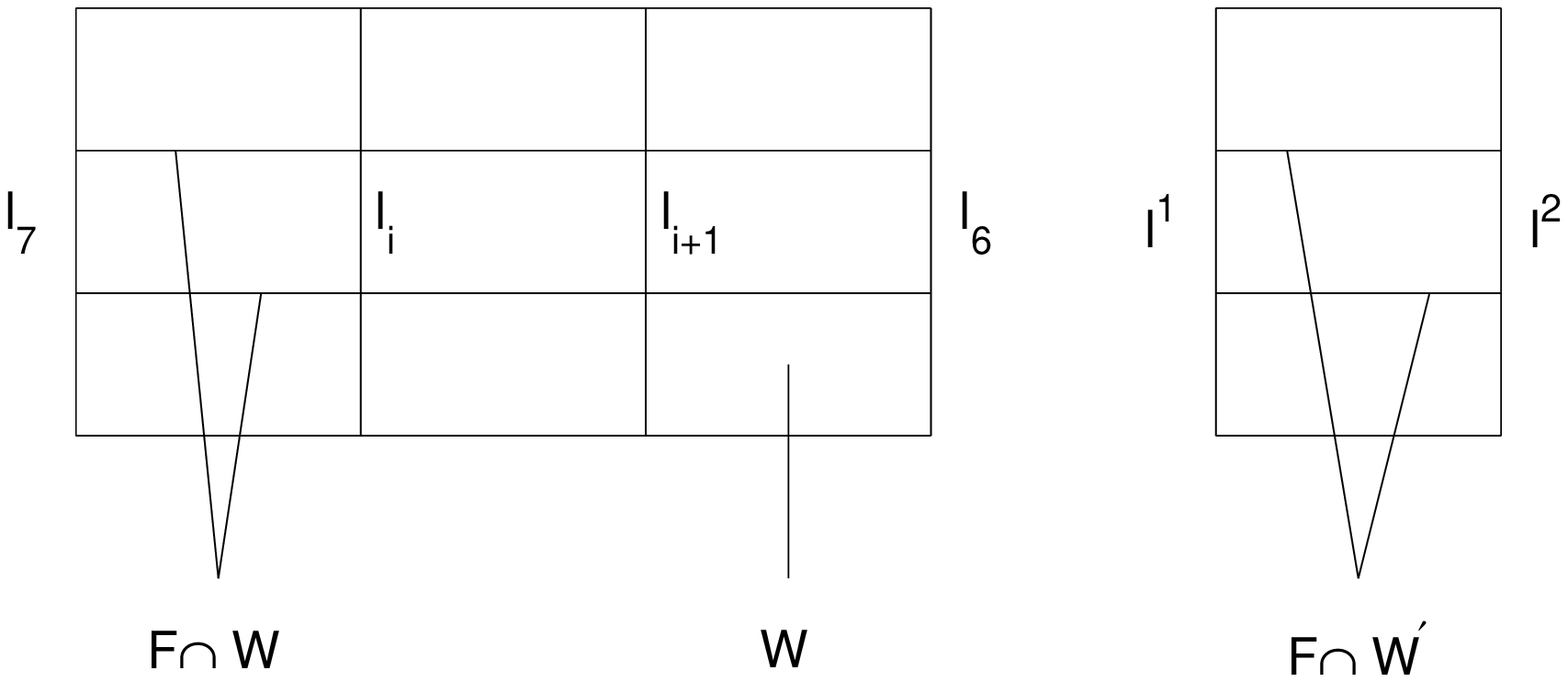}
\begin{center}
Figure 5.1
\end{center}
\end{center}

{\bf Lemma 5.0.} \ (1) Each component of $F\cap(F_{1}\cup F_{2})$
is an essential circle in both $F$ and $F_{1}\cup F_{2}$.

(2) Each component of  $F\cap W$ is an arc in $W$ such that one of
its two endpoints lies in $l_{6}$ and the other lies in $l_{7}$.
Similarly each component of $F\cap W^{'}$ is an arc in $W^{'}$
such that one of its two endpoints lies in $l^{1}$ and the other
lies in $l^{2}$. Hence $|F\cap l_{i}|=|F\cap l_{j}|$ for all $i,j$
and $|F\cap l^{1}|=|F\cap l^{2}|$ as in Figure 5.1.

(3) Each component of $F\cap (F_{1}\cup F_{2})$ isotopic to
$\partial A_{i}$ is disjoint from $W\cup W^{'}$.

For two surfaces $P_{1}$ and $P_{2}$ in a 3-manifold, a pattern of
$P_{1}\cap P_2$ is a set of disjoint arcs and circles representing
isotopy classes of $P_{1}\cap P_{2}$. For each isotopy class $s$,
we use $\nu(s)$ to denote the number of components of $P_{1}\cap
P_2$ in the isotopy class $s$.
\begin{center}
\includegraphics[totalheight=5cm]{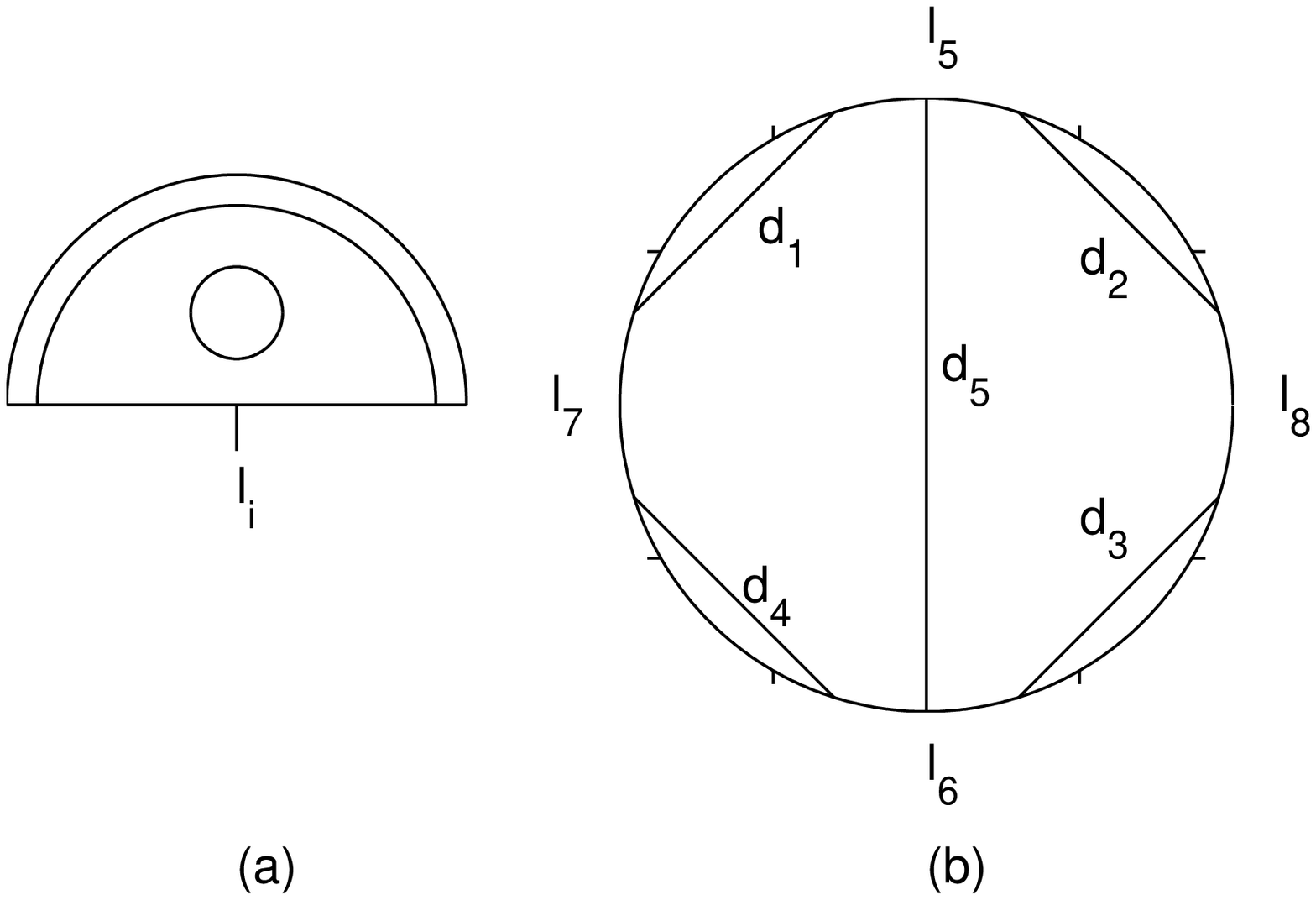}
\begin{center}
Figure 5.2
\end{center}
\end{center}
\begin{center}
\includegraphics[totalheight=5cm]{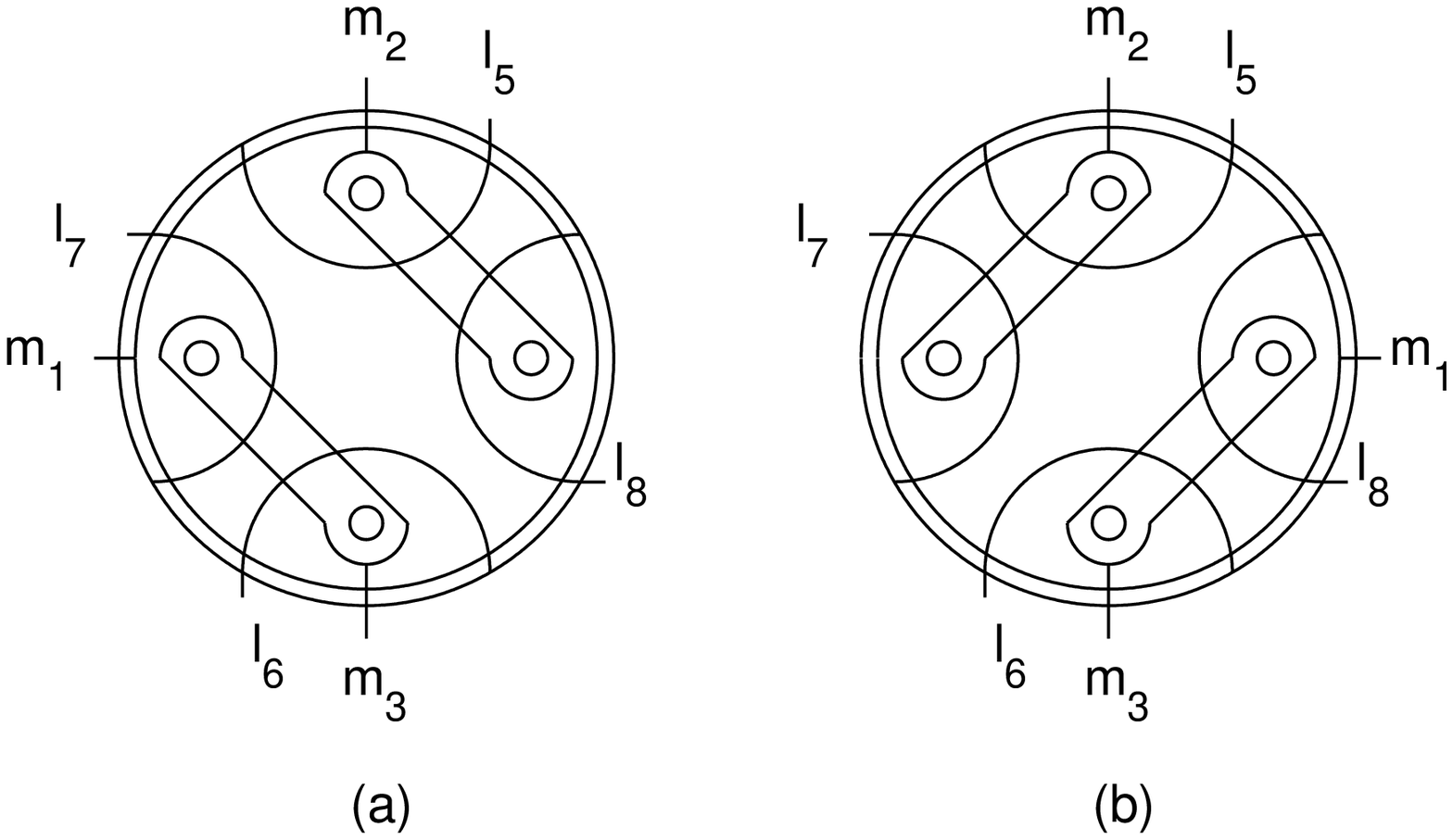}
\begin{center}
Figure 5.3
\end{center}
\end{center}

Lemma 5.1 follows also immediately from the proof of Lemma 4.3 in
[QW2].

{\bf Lemma  5.1.} \ Each component of $F\cap M_{3}$ is isotopic to
one of $\partial H\cap M_{3}$, $A_{5}$ and $A_{7}$.

{\bf Proof.} \ Note that the four arcs $l_{5}, l_{6}, l_{7},
l_{8}$ separate $F_{2}$ into four annulus $A^{5}, A^{6}, A^{7},
A^{8}$ and a disk $D$. By the minimality of $|F\cap W|$, the
pattern of $F\cap A^{j}$ is as in Figure 5.2(a) and the pattern of
$F\cap D$ is as in Figure 5.2(b). Since $|F\cap l_{i}|$ is a
constant, $\nu(d_{5})=0$. If $\nu(d_{i})\neq 0$ for $1\leq i\leq
4$, then $F\cap F_{2}$ contains $min\bigl\{\nu
(d_{1}),\ldots,\nu(d_{4})\bigr\}$ components parallel to a disk on
$\partial E_{2}$. Now if  $\nu (d_{1})=0$, then $\nu(d_{3})=0$.
Similarly, if $\nu(d_{2})=0$, then $\nu(d_{4})=0$. Thus according
to the order of $l_{5},l_{6},l_{7},l_{8}$ in $F_{2}$, the pattern
of $F\cap F_{2}$ is as in one of Figure 5.3(a) and (b) with
$\nu(m_{2})=\nu(m_{3})$. Note that $W_{5}$ and $W_{7}$ separate
$M_{3}$ into three solid tori $J^{1}, J^{2}, J^{3}$. Without loss
of generality, we assume that $A_{5}\subset J^{1}, A_{7}\subset
J^{2}$. Let $S=F\cap M_{3}$ and $S^{'}$ be a component of $S$.

\begin{center}
\includegraphics[totalheight=5cm]{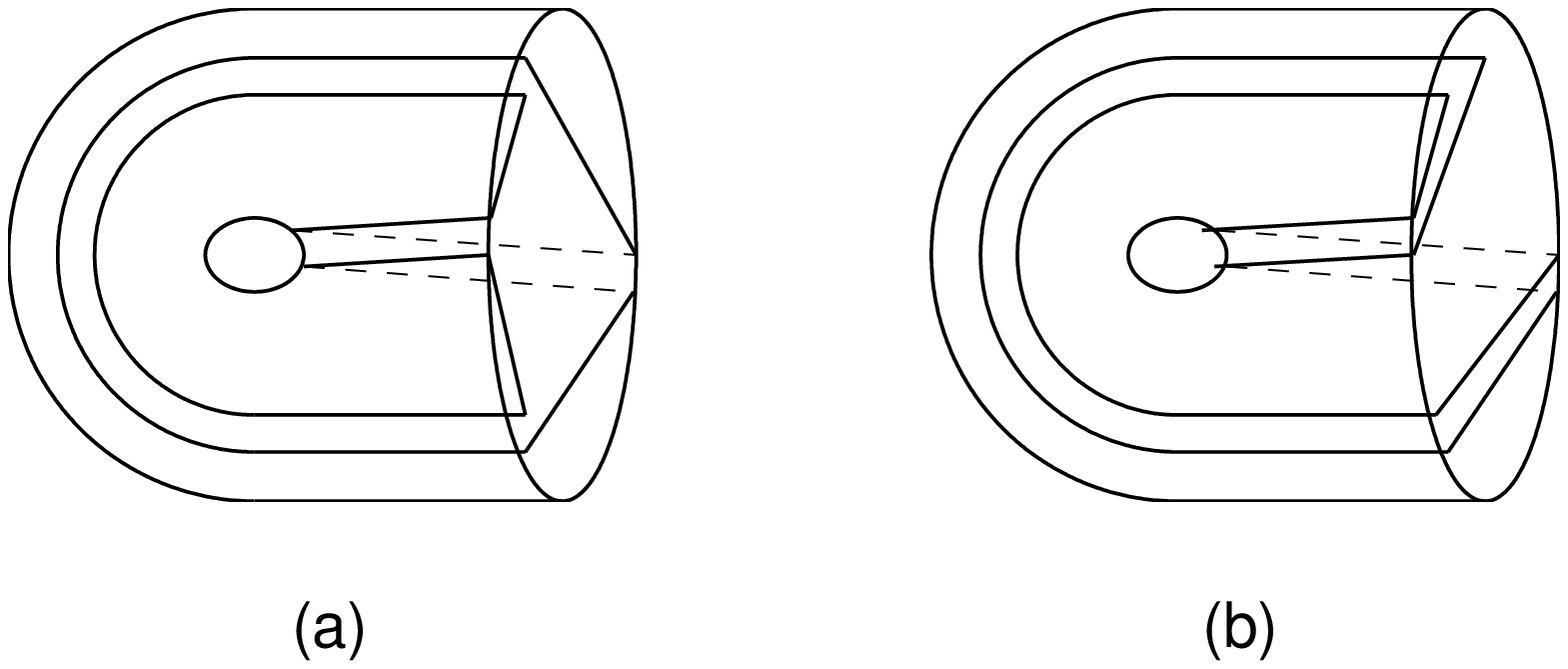}
\begin{center}
Figure 5.4
\end{center}
\end{center}

Now we claim that if one of component of $\partial S^{'}$ is
isotopic to $\partial E_{2}$, then $S^{'}$ is isotopic to
$\partial H\cap M_{3}$.

Let $\partial_{1} S$ be the outermost component of $\partial S$
isotopic to $\partial E_{2}$. Now $\partial _{1} S$ intersects
$l_{i}$ as in Figure 5.3(a) or (b). Without loss of generality, we
assume that $\partial_{1} S\subset \partial S^{'}$. We denote by
$e_{i}$ the arc $\partial_{1} S\cap A^{i}$. Now let $S_l=S^{'}\cap
J^l$, then $S_l$ is an incompressible surface in $J^{l}$. Note
that $\partial S_1=e_{5}\cup e_{6}\cup (S\cap W_{5})$ bounds a
disk in $J^{1}$ parallel to a disk on $\partial M_{3}$. Similarly
$S_2$ is a disk in $J^{2}$ parallel to a disk on $\partial M_{3}$
bounded by $e_{7}\cup e_{8}\cup (S\cap W_{7})$. $\partial S_3$
also has one component which is trivial in $\partial M_{3}$ as in
Figure 5.4(a). Hence one component of $S_3$ is a disk in $J^{3}$
parallel to $\partial J^{3}$, Thus $S^{'}=S_1\cup_{S\cap W_{5}}
S_3\cup _{S\cap W_{7}} S_2$ is isotopic to $M_{3}\cap\partial H$.

Now we claim that $\nu(m_2)=\nu(m_3)=0$ in Figure 5.3(a) and (b).

Let $S_{0}=S-X^{'}$ where $X^{'}$ is a subset of $S$ such that
each component of $X^{'}$ is isotopic to $\partial H\cap M_{3}$.
Then each component of $\partial S_{0}$ is not isotopic to
$\partial E_{2}$.  Let $P_{3}=S_{0}\cap J^{3}$. If $\nu(m_{2})\neq
0$, Then $P_{3}$ is incompressible in $J^{3}$ and $\partial P_{3}$
contains $2\nu(m_{2})=2\nu(m_{3})$ components $c$ as in Figure
5.4(b). Since $a_{7}$ intersects a basis disk $B_{3}$ of $J_{3}$
in three points and $a_{5}$ intersects $B_{3}$ in one points, $c$
does not bound a disk in $J^{3}$.  Since $J^{3}$ is a solid torus,
each component of $P_{3}$ is an annulus which is
$\partial$-compressible. Let $D^*$ be a $\partial$-compressing
disk of an outermost  component of $P_{3}$. Note that the
$\partial$-compressing disk $D^*$ can be isotoped so that
$D^*\cap\partial J^3\subset E_2\cap J^3$. Then back to $J_3$,
$D^*$ is isotopic to one of $D^{1}, D^{2}, D^{3}$ as in Figure
5.5. In the case of $D^{1}$ or $D^{2}$, one can push $F$ along the
disc to reduce $|F\cap W|$; In the case of $D^{3}$, one can push
$F$ along the disc to reduce $|F\cap  F_{2}|$, but not to increase
$|F\cap W|$. In each case, it contradicts the minimality of
$C(F)$.
\begin{center}
\includegraphics[totalheight=6cm]{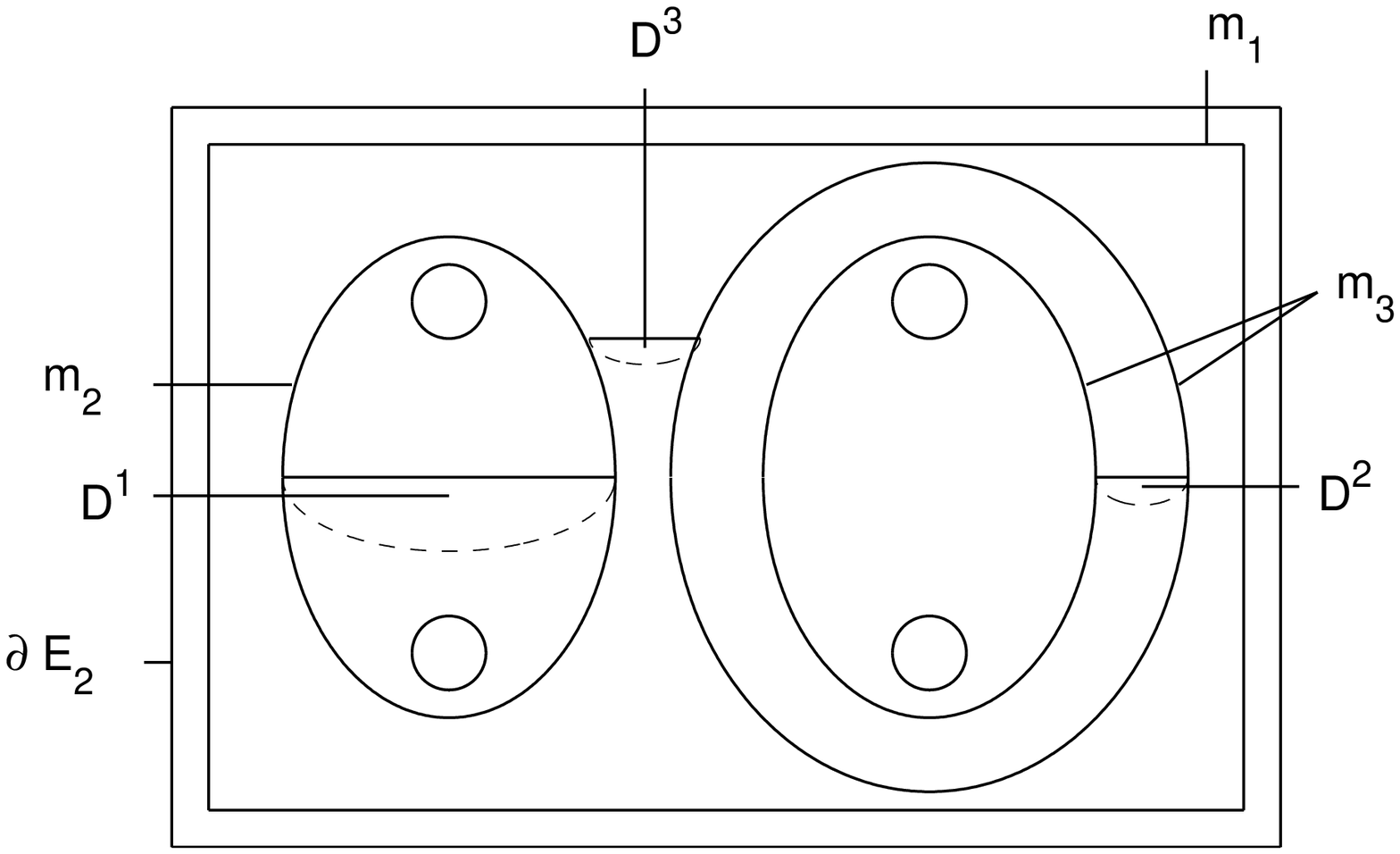}
\begin{center}
Figure 5.5
\end{center}
\end{center}

Now let $P$ be a component of $S=F\cap M_{3}$. If one component of
$\partial P$ is isotopic to $\partial E_{2}$, then $P$ is isotopic
to $M_{3}\cap \partial H$. If not, then each component of
$\partial P$ is isotopic to one component of $\partial A_{5}\cup
\partial A_{7}$. By the minimality of $C(F)$, $P$ is contained
in $J^{1}$ or $J^{2}$. It is easy to see that $P$ is isotopic to
one of $A_{5}$ and $A_{7}$.\qquad Q.E.D.

Now we consider $S=F\cap M_{1}$. Note that $W_{1}, W^{'}$ separate
$M_{1}$ into two solid tori $J^{1}, J^{2}$ and a handlebody of
genus two $H^{'}$ such that $A_{1}\subset J^{1}$,
$A_{3},A_{9}\subset H^{'}$; $l_{1},l_{2},l^{1},l^{2}$ separate
$F_{1}$ into two annuli and two planar surfaces with three
boundary components and a disk $D$ such that $\partial J^{2}\cap
F_{1}=D$, see Figure 5.6. Let $k_{1}(k_{2})$ be a component of
$F\cap W_{1} (F\cap W^{'})$, and $k_{1}^{'}(k_{2}^{'})$ be an arc
in $D$ connecting the two endpoints of $k_{1}(k_{2})$. Let
$\alpha=k_{1}\cup k_{1}^{'}$ and $\beta=k_{2}\cup k_{2}^{'}$. Note
that  $k_{1}^{'}$ and $k_{2}^{'}$ can be chosen so that $\beta$
intersects $\alpha$ in one point. Furthermore, by construction,
$\alpha$ intersects a basis disk of $J^{2}$ in two points and
$\beta$ intersects a basis disk of $J^{2}$ in one point. Now  we
fix the orientations of $\alpha$ and $\beta$ so that
$\alpha=y^{2}$ and $\beta=y$ where $y$ is a generator of
$\pi_{1}(J^{2})$. Then $\alpha\beta^{-2}$ is an essential circle
in $\partial J^{2}$ which  is null homotopic  in $J^{2}$.

{\bf Lemma 5.2.} \ Let $P$ be a component of $S=F\cap M_{1}$. If
one component of $\partial P$ is isotopic to $\partial E_{1}$.
Then $P$ is isotopic to $M_{1}\cap
\partial H$.

{\bf Proof.} \ This follows immediately from the proof of Lemma
5.1.\qquad Q.E.D.

By the construction and Lemma 5.0, the pattern of $\partial S\cap
(F_{1}\cap (J^{1}\cup H^{'}))$ is as in one of Figure 5.6(a) and
(b) such that

1) in Figure 5.6(a), $\nu(f_{1})=\nu(f_{2}),
\nu(f_{3})=\nu(f_{5}), \nu(f_{4})=\nu(f_{6})$ and
$\nu(f_{3})+\nu(f_{4})=\nu(f_{1})$.

2) in Figure 5.6(b),
$\nu(f_{1})=\nu(f_{2})=\nu(f_{3})+\nu(f_{4})$,
$\nu(f_{3})=\nu(f_{6})$ and
$\nu(f_{4})=\nu(f_{5})=\nu(f_{7})=\nu(f_{8})\neq 0$.

\begin{center}
\includegraphics[totalheight=4.5cm]{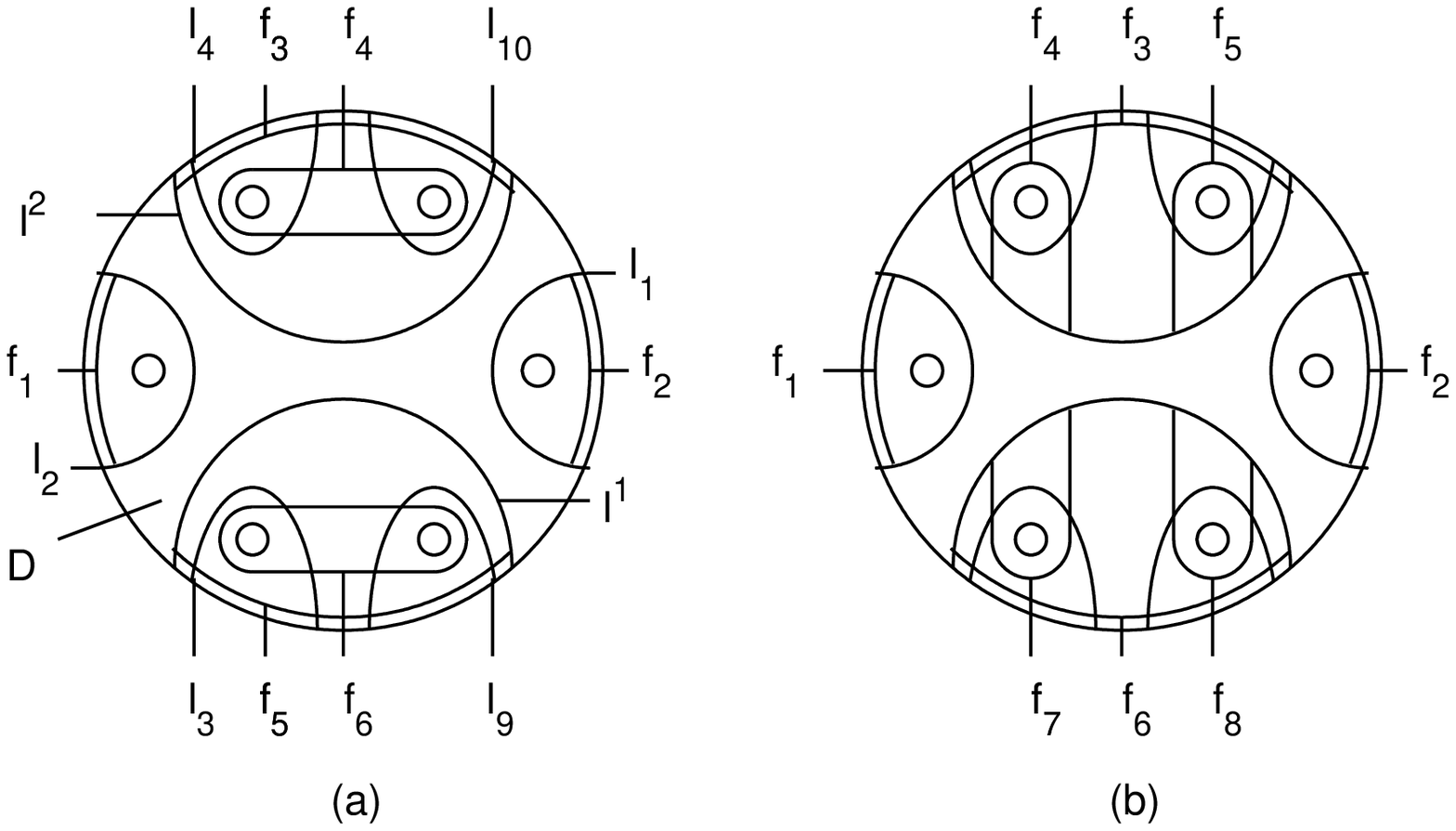}
\begin{center}
Figure 5.6
\end{center}
\end{center}
\begin{center}
\includegraphics[totalheight=5cm]{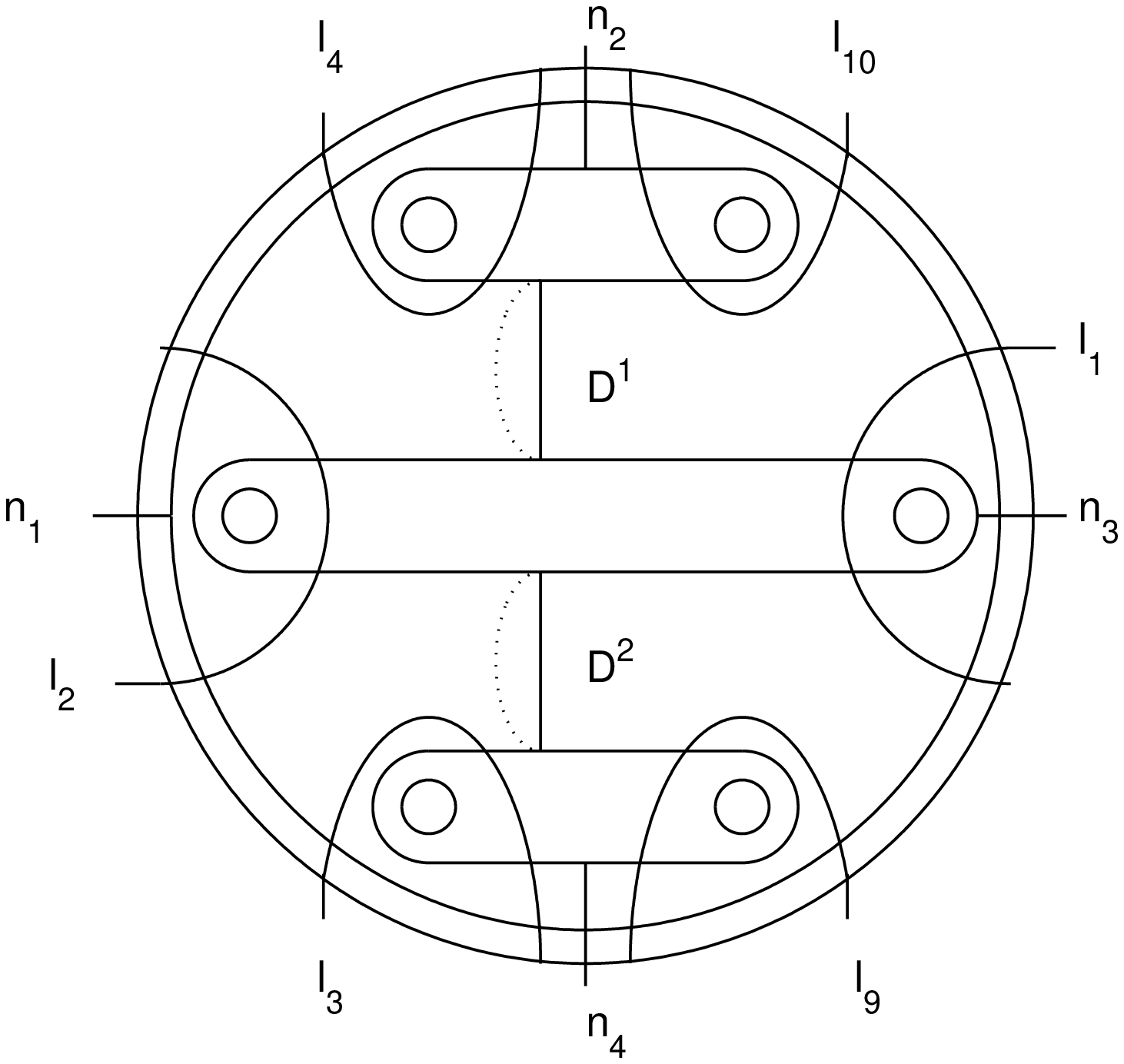}
\begin{center}
Figure 5.7
\end{center}
\end{center}

{\bf Lemma 5.3.} \ If the pattern of $S\cap (F_{1}\cap (J^{1}\cup
H^{'}))$ is as in Figure 5.6(a), then  the pattern of $S\cap
F_{1}$ is as in Figure 5.7 with
$\nu(n_{2})=\nu(n_{3})=\nu(n_{4})$.

{\bf Proof.} \ If $\nu(f_{3})=0$ in Figure 5.6(a), then the
pattern of $S\cap F_{1}$ is as in Figure 5.7 with
$\nu(n_{2})=\nu(n_{3})=\nu(n_{4})$ and $\nu(n_{1})=0$.

Suppose now that $\nu(f_{3})\neq 0$ in Figure 5.6(a).  Since
$\nu(f_{3})=\nu(f_{5})\leq \nu(f_{1})=\nu(f_{2})$, the pattern of
$S\cap D$ is as in Figure 5.8 where $\nu(d_{1})=\nu(d_{3})$ and
$\nu(d_{2})=\nu(d_{4})$.  If $\nu(d_{1}), \nu(d_{2})\neq 0$, then
$S\cap F_{1}$ contains $min(\nu(d_{1}),\nu(d_{2}))$ components
isotopic to $\partial E_{1}$. Thus if $\nu(d_{1})=\nu(d_{2})$,
then $S\cap F_{1}$ is as in Figure 5.7 with
$\nu(n_{2})=\nu(n_{3})=\nu(n_{4})$.  Now without loss of
generality, we assume that $\nu(d_{1})<\nu(d_{2})$. Let
$k=\nu(d_{2})-\nu(d_{1})$. By Lemma 5.0(2) and Lemma 5.2,
$\partial (S\cap J^{2})$ contains $n=gcd(k,k+\nu(d_{5}))$
components $c$ isotopic to $\alpha^{p}\beta^{q}$ where
$|p|=(k+\nu(d_{5}))/n$ and $|q|=k/n$. Since $y+\nu(d_{5})\geq y$,
$c$ is not null homotopic in $J^{2}$. Furthermore, $c\cap
d_{2}\neq \phi$, $c\cap d_{4}\neq \phi$; if $\nu(d_{5})\neq 0$,
then $c\cap d_{5}\neq \phi$. Thus these curves separates $\partial
J^{2}$ into $m$ annuli $A^{1},\ldots,A^{m}$ such that, for each
$j$, there is an arc in $D\cap A^{j}$ connecting the two boundary
components of $A^{j}$. Since $J^{2}$ is a solid torus, each
component of $(S-X(F))\cap J^{2}$ is an annulus. Let $D^{*}$ be a
$\partial$-compressing disk of $(S-X(F))\cap J^{2}$. Then $D^{*}$
can be moved so that $D^{*}\cap
\partial J^{2}=D^{*}\cap D=a$. Thus there are three possibilities:
\begin{center}
\includegraphics[totalheight=6cm]{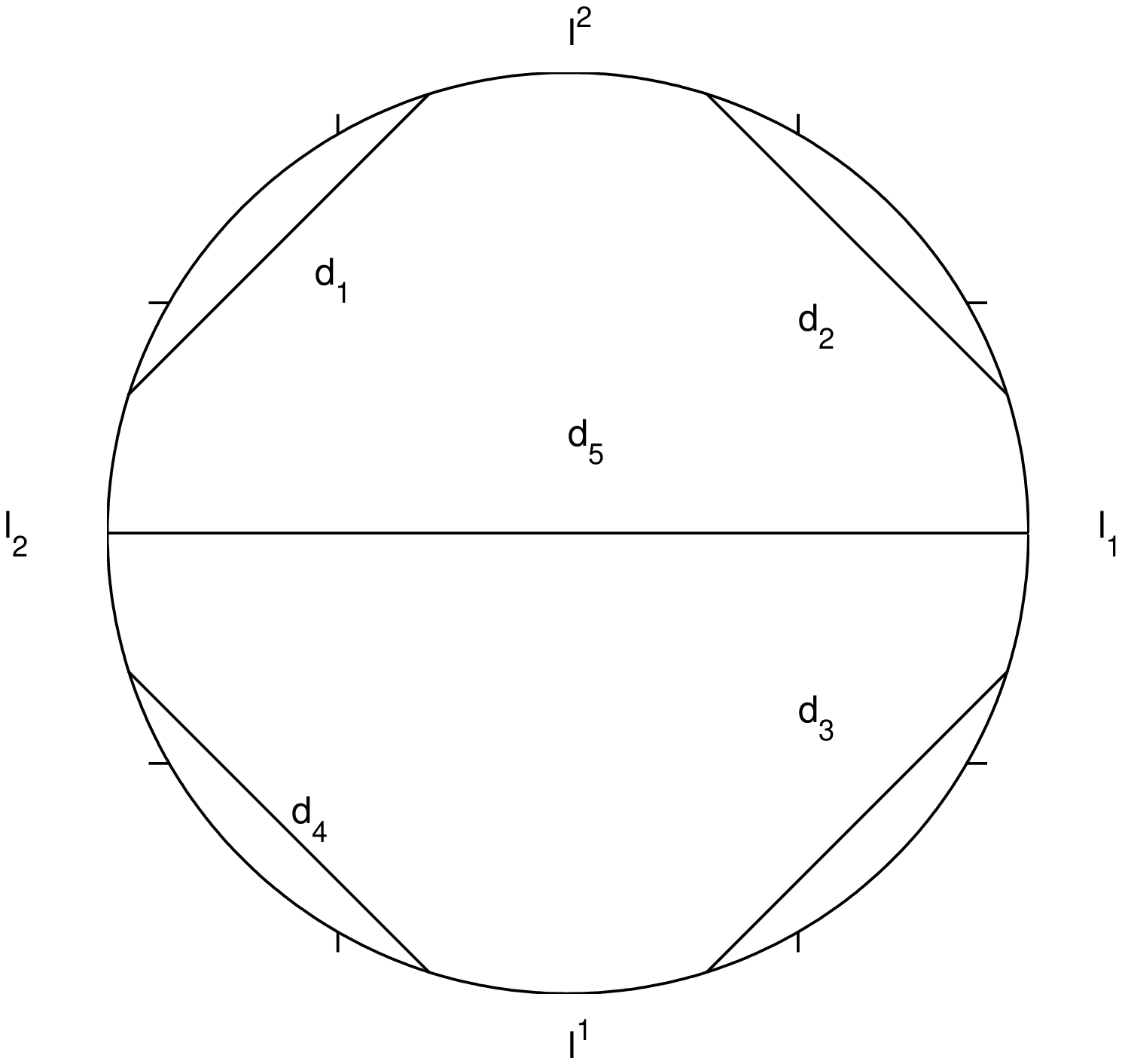}
\begin{center}
Figure 5.8
\end{center}
\end{center}

1) the two endpoints of $a$ lie in one of $d_{2},d_{4},d_{5}$. Now
$D^{*}$ is one of $D^{1},D^{3}$ as in Figure 5.9. In each case,
one can push $F$ along the disc to reduce $|F\cap W|$, a
contradiction.

2) one endpoint of $a$ lies in $d_{2}\cup d_{4}$ and the other
lies in $d_{5}$. Now $D^{*}$ is  $D^{2}$ as in Figure 5.9. This
case is similar to Case 1).

3) one endpoint of $a$ lies in $d_{2}$ and the other lies in
$d_{4}$.  In this case, $\nu(d_{5})=0$. By Lemma 5.0(2),
$\nu(f_{4})=\nu(f_{6})=0$ in Figure 5.6(a).  Now the pattern of
$S\cap F_{1}$ is as in Figure 5.10 and  $D^{*}$ is as  in Figure
5.10. By doing a surgery on $F$ along $D^{*}$, we can obtain a
surface $F^{'}$ isotopic to $F$ such that $|F^{'}\cap W|=|F\cap
W|$, $|F^{'}\cap F_{2}|=|F\cap F_{2}|$, $|(F^{'}\cap
M_{1}-X(F^{'}))\cap W^{'}|<|(F\cap M_{1}-X(F))\cap W^{'}|$ (by
Lemma 5.2), a contradiction.\qquad Q.E.D.

\begin{center}
\includegraphics[totalheight=7cm]{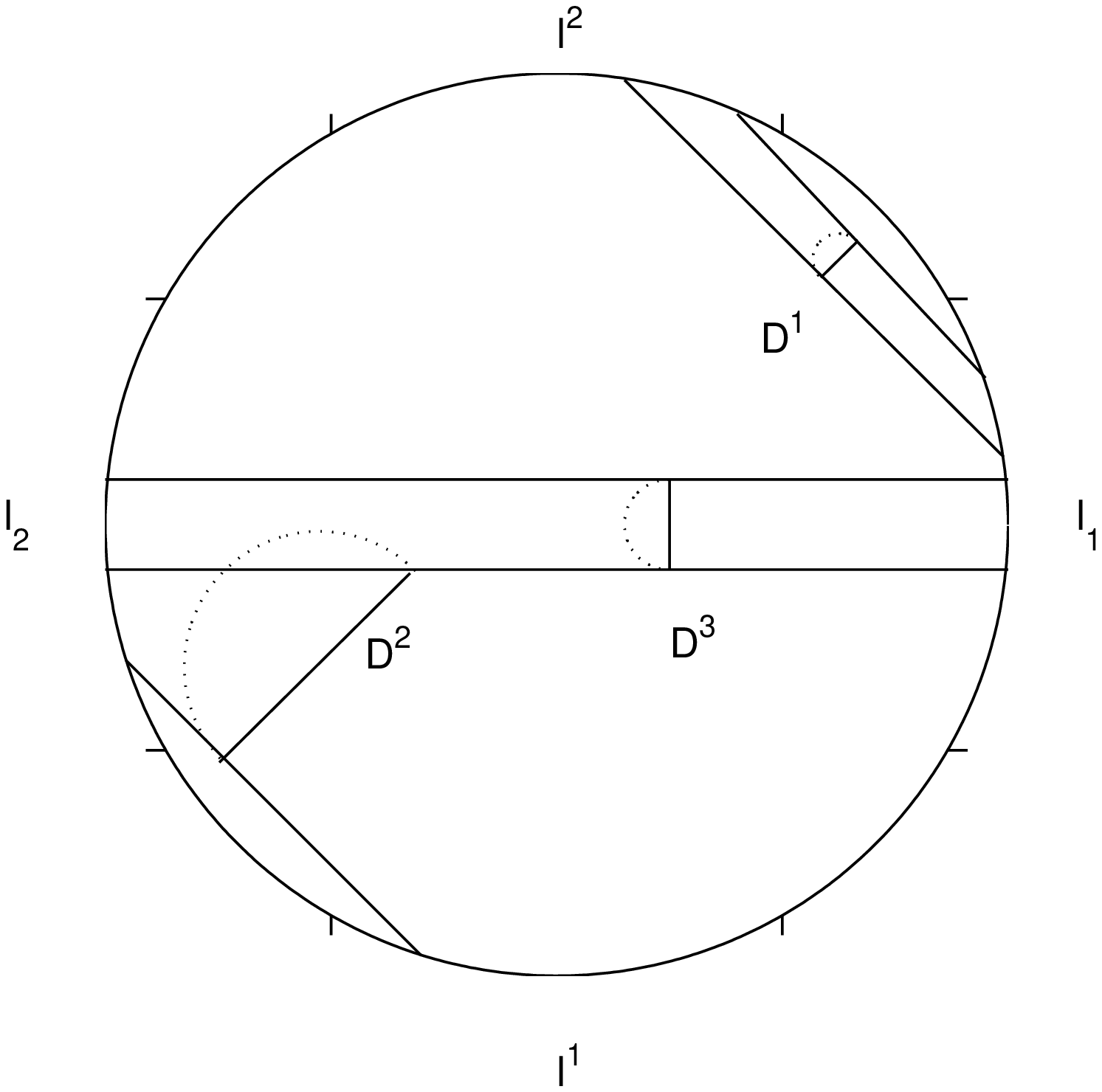}
\begin{center}
Figure 5.9
\end{center}
\end{center}

\begin{center}
\includegraphics[totalheight=5cm]{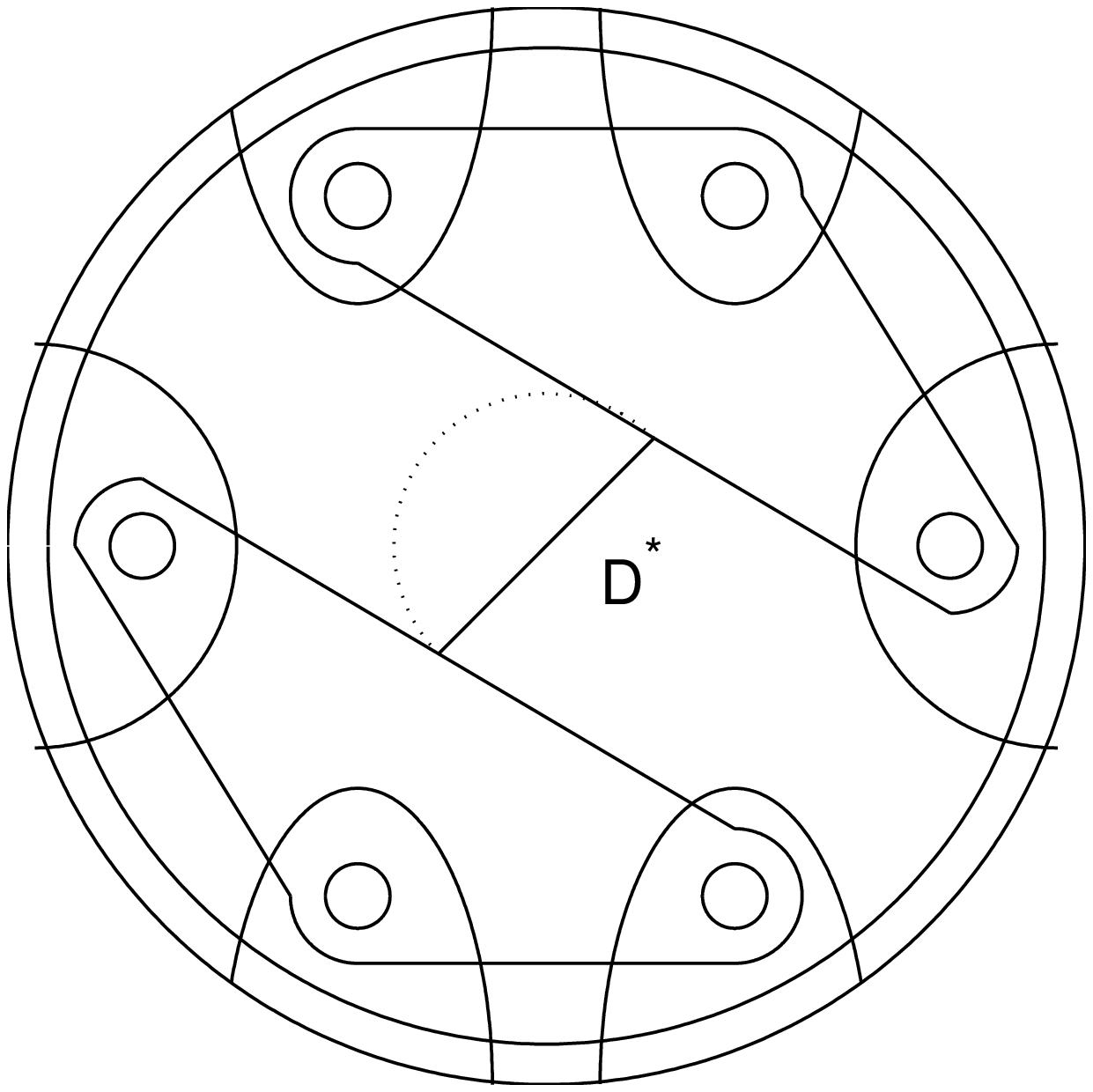}
\begin{center}
Figure 5.10
\end{center}
\end{center}

{\bf Lemma 5.4.} \ If the pattern of $S\cap (F_{1}\cap (J_{1}\cup
H^{'}))$ is as in Figure 5.6(b), then the pattern of $S\cap F_{1}$
is as in Figure 5.11.

{\bf Proof.} \ In this case,
$\nu(f_{1})=\nu(f_{2})=\nu(f_{3})+\nu(f_{4})=\nu(f_{6})+\nu(f_{7})$.
Thus the pattern of $S\cap D$ is as in Figure 5.12 where
$\nu(d_{1})=\nu(d_{3})$, $\nu(d_{2})=\nu(d_{4})$, and
$\nu(d_{5})=2\nu(f_{5})$.
\begin{center}
\includegraphics[totalheight=6.5cm]{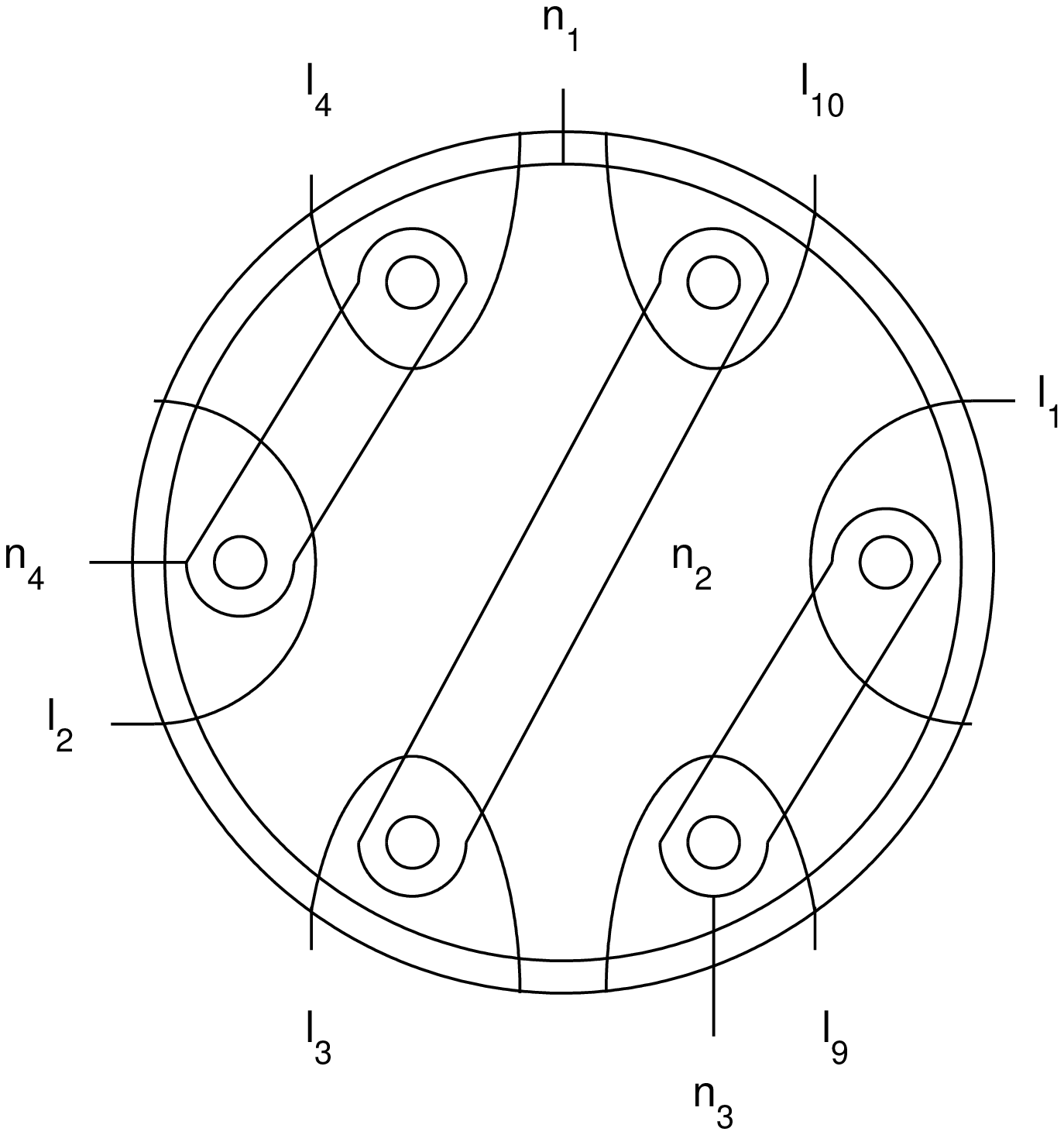}

\begin{center}
Figure 5.11
\end{center}
\end{center}
\begin{center}
\includegraphics[totalheight=6cm]{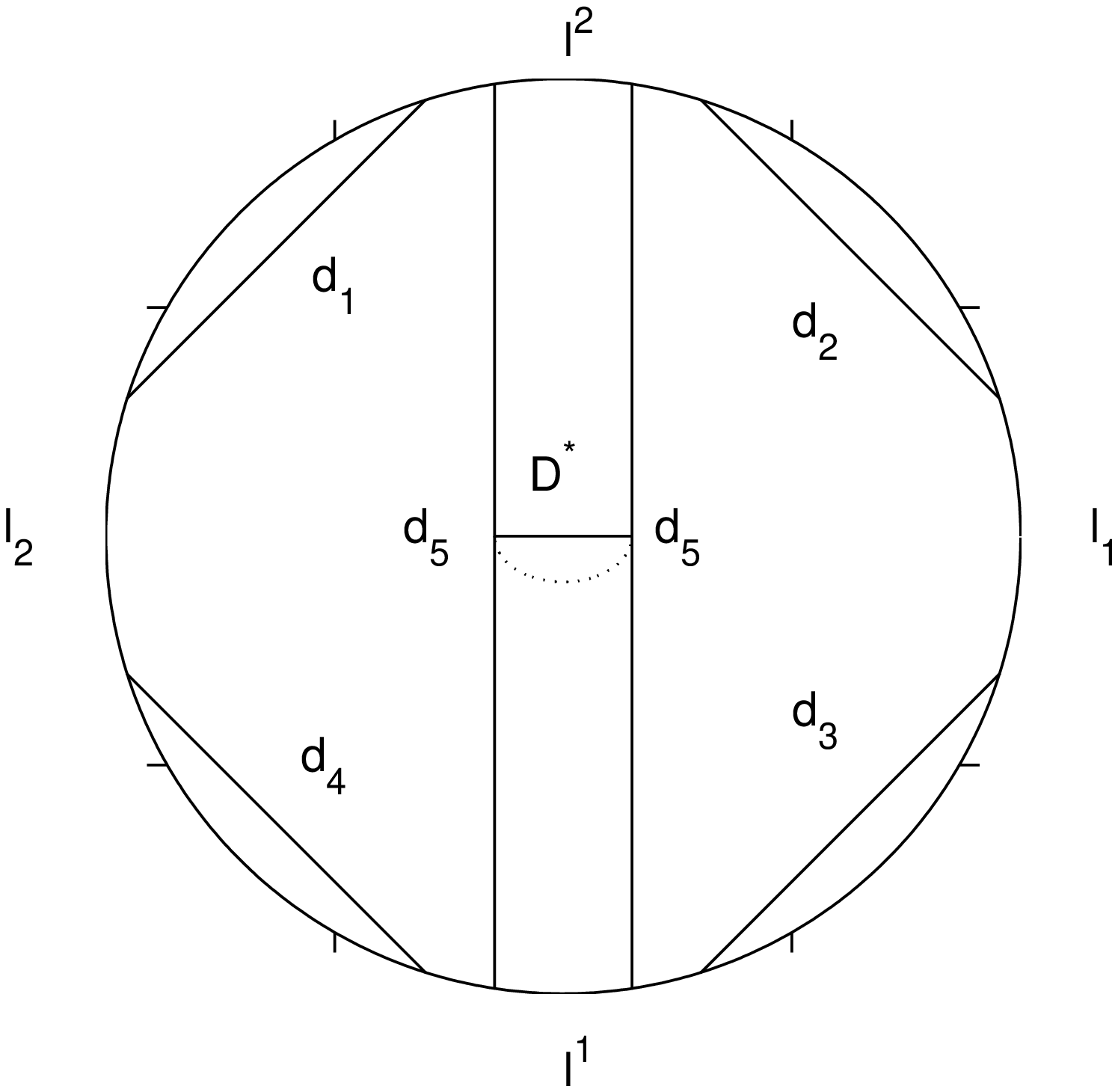}
\begin{center}
Figure 5.12
\end{center}
\end{center}

Now $\nu(d_{5})\neq 0$. There are two cases:

1. \ $\nu(f_{3})=\nu(f_{6})= 0$ in Figure 5.6(b).

In this case, $\nu(d_{5})=\nu(d_{1})+\nu(d_{2})$. Now there are
two subcases:

1). $\nu(d_{1})=\nu(d_{2})$.

Since $\nu(d_{5})\neq 0$,  $\partial (S\cap J^{2})$ contains
$\nu(d_{1})$ trivial components in $\partial J^{2}$ which bound
some disks in $S$ as in Figure 5.4(a) and $\nu(d_{5})$ components
isotopic to $\beta$. Since $\beta$ intersects a basis disk of
$J^{2}$ in one point, each non-trivial component of $S\cap J^{2}$,
say $A^{*}$, is an annulus which is parallel to each component of
$\partial J^{2}-\partial A^{*}$. Thus there is a
$\partial$-compressing disk $D^{*}$ of $S\cap J^{2}$ as in Figure
5.12. Thus by doing a surgery on $F$ along $D^{*}$, we can obtain
a surface $F^{'}$ isotopic to $F$ such that $|F^{'}\cap W|=|F\cap
W|$, $|F^{'}\cap F_{2}|=|F\cap F_{2}|$, $|(F^{'}\cap
M_{1}-X(F^{'}))\cap W^{'}|<|(F\cap M_{1} -X(F))\cap W^{'}|$, a
contradiction.

2)  $\nu(d_{1})\neq\nu(d_{2})$.

Let $k=|\nu(d_{2})-\nu(d_{1}|$ and $n=gcd(k,k+\nu(d_{5})$.

We first suppose  that $\nu(d_{1})<\nu(d_{2})$.  Then $\partial
(S\cap J^{2})$ contains $\nu(d_{1})$ trivial components and $n$
components $c$ isotopic to $\alpha^{p}\beta^{q}$ where
$|q|=(k+\nu(d_{5}))/n$  and $|p|=k/n$. By the construction, $p>0$
if and only if $q>0$. (See Figure 2.2.) That means that $c$ is not
null homotopic in $J^{2}$. By the proof of Lemma 5.3, we can
obtain a surface $F^{'}$ isotopic to $F$ such that
$C(F^{'})<C(F)$, a contradiction.

Suppose that $\nu(d_{1})>\nu(d_{2})$.  By the above argument,
$\partial (S\cap J^{2})$ contains $\nu(d_{2})$ trivial components
and $n$ components $c$ isotopic to $\alpha^{p}\beta^{q}$ where
$|q|=(k+\nu(d_{5}))/n$ and $|p|=k/n$. If $c$ is not null homotopic
in $J_{2}$, then by the above argument, we can obtain a surface
$F^{'}$ isotopic to $F$ so that $C(F^{'})<C(F)$, a contradiction.
Assume that $q=-2p$. Then $\nu(d_{5})=\nu(d_{1})-\nu(d_{2})$.
Since $\nu(d_{5})=\nu(d_{1})+\nu(d_{2})$, $\nu(d_{2})=0$ and
$\nu(d_{5})=\nu(d_{1})$. Thus $F_{1}\cap F$ is as in Figure 5.11
with $\nu(n_{2})=\nu(n_{3})=\nu(n_{4})$ and $\nu(n_{1})=0$.

2. \ $\nu(f_{3})=\nu(f_{6})\neq 0$ in Figure 5.6(b).

Now there are two subcases:

1) $\nu(d_{1})\leq \nu(d_{2})$.

Now $S\cap F_{1}$ contains $min(\nu(d_{1}),\nu(f_{3}))$ components
isotopic to $\partial E_{1}$. If $\nu(d_{1})\geq \nu(f_{3})$, then
by the argument in Case 1, we can obtain a surface $F^{'}$
isotopic to $F$ such that $C(F^{'})<C(F)$, a contradiction. Assume
that $\nu(d_{1})<\nu(f_{3})$, then $S\cap F_{1}$ contains
$\nu(d_{1})$ components isotopic to $\partial E_{1}$. Now
$2\nu(f_{1})=\nu(d_{1})+\nu(d_{2})$. By assumption,
$\nu(f_{1})=\nu(f_{3})+\nu(f_{4})$. Thus $\nu(d_{1})<\nu(d_{2})$.
Then, by the proof of Lemma 5.3, $\partial (S\cap J^{2})$ contains
$gcd(k,k+\nu(d_{5}))$ components each of which is isotopic to
$\alpha^{p}\beta^{q}$ where $|q|=(k+\nu(d_{5}))/n$ and $|p|=k/n$,
where  $k=|\nu(d_{2})-\nu(d_{1}|$ and $n=gcd(k,k+\nu(d_{5})$.  If
$q\neq -2p$, then by the proof of Lemma 5.3, in $H_{K}$, there is
an essential closed surface $F^{'}$ isotopic to $F$ such that
$C(F^{'})<C(F)$, a contradiction. Since
$y=\nu(d_{2})-\nu(d_{1})=2(\nu(f_{1})-\nu(d_{1}))>
2(\nu(f_{1})-\nu(f_{3}))=2\nu(f_{5})$, $\nu(d_{5})=2\nu(f_{5})$.
Thus $q\neq -2p$.

2) $\nu(d_{1})> \nu(d_{2})$.

Now $S\cap F_{1}$ contains $min(\nu(d_{2}),\nu(f_{3}))$ components
isotopic to $\partial E_{1}$. If $\nu(d_{2})\geq \nu(f_{3})$, then
by the argument in Case 1, the pattern of $F\cap F_{1}$ is as in
Figure 5.11 with $\nu(n_{1})=\nu(f_{3})$ and
$\nu(n_{2})=\nu(n_{3})=\nu(n_{4})$. By the above argument, it is
impossible that $\nu(d_{1})<\nu(f_{3})$.\qquad Q.E.D.
\begin{center}
\includegraphics[totalheight=5.5cm]{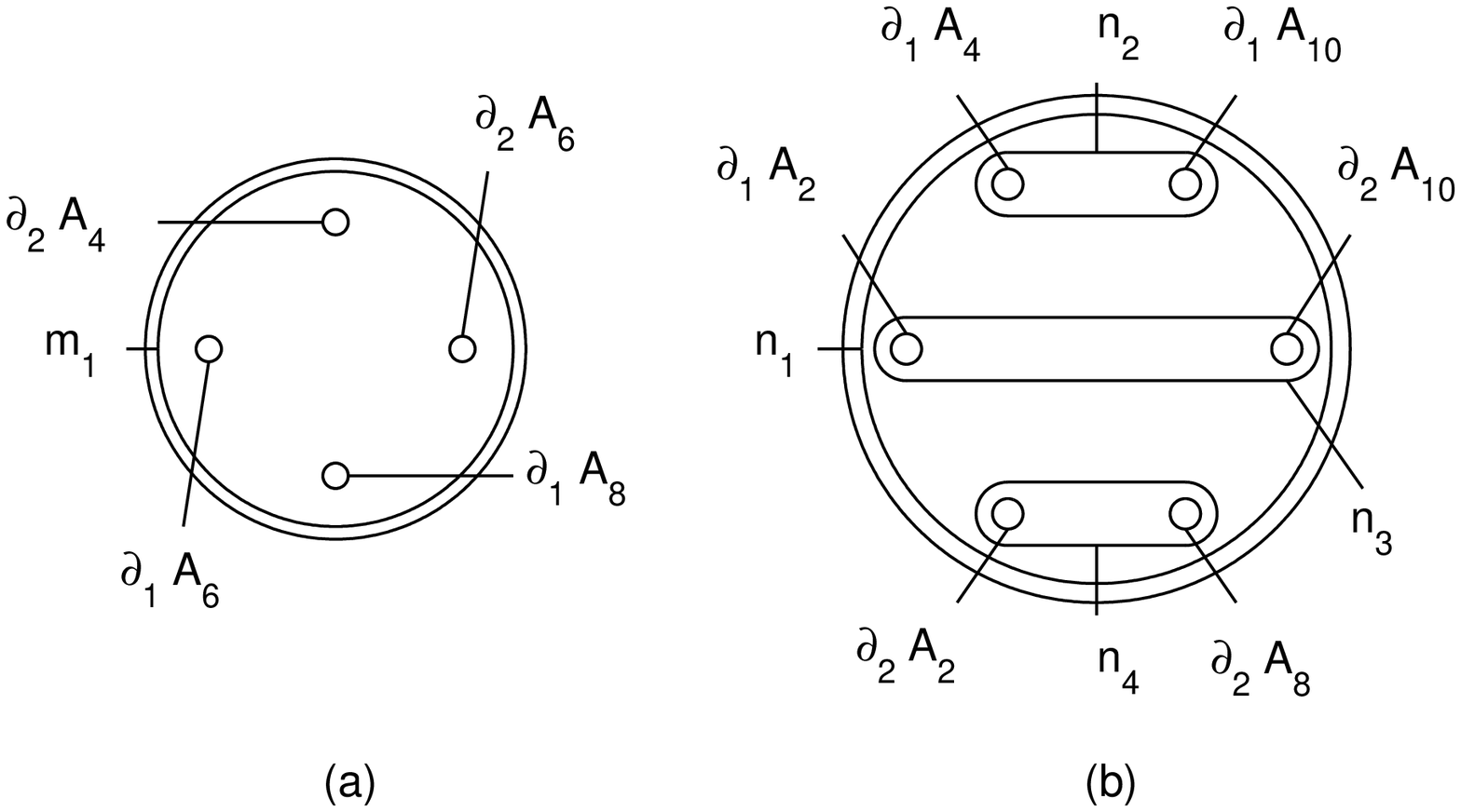}
\begin{center}
Figure 5.13
\end{center}
\end{center}
\begin{center}
\includegraphics[totalheight=5.5cm]{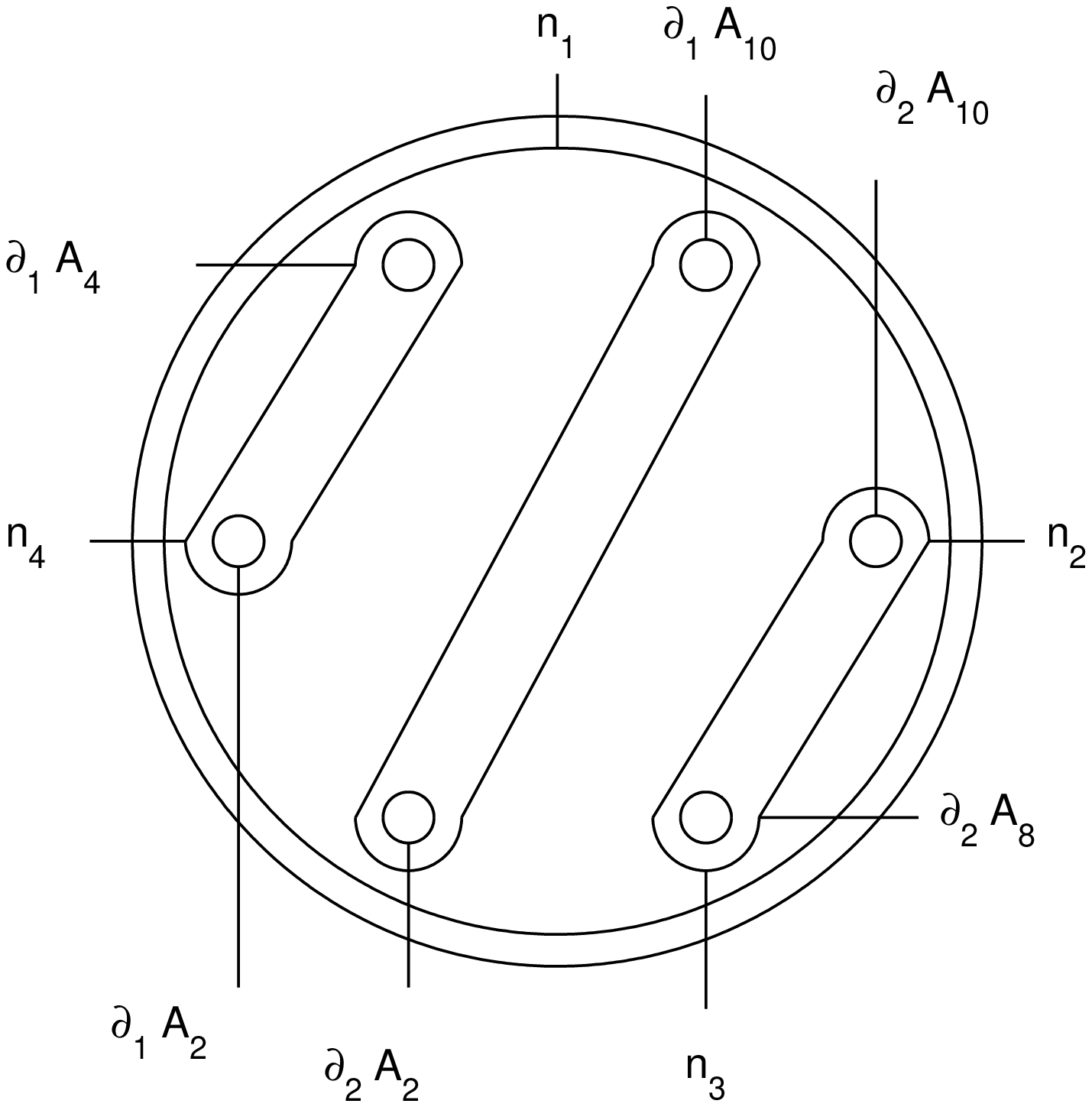}
\begin{center}
(c)
\end{center}
\begin{center}
Figure 5.13
\end{center}
\end{center}
\begin{center}
\includegraphics[totalheight=5cm]{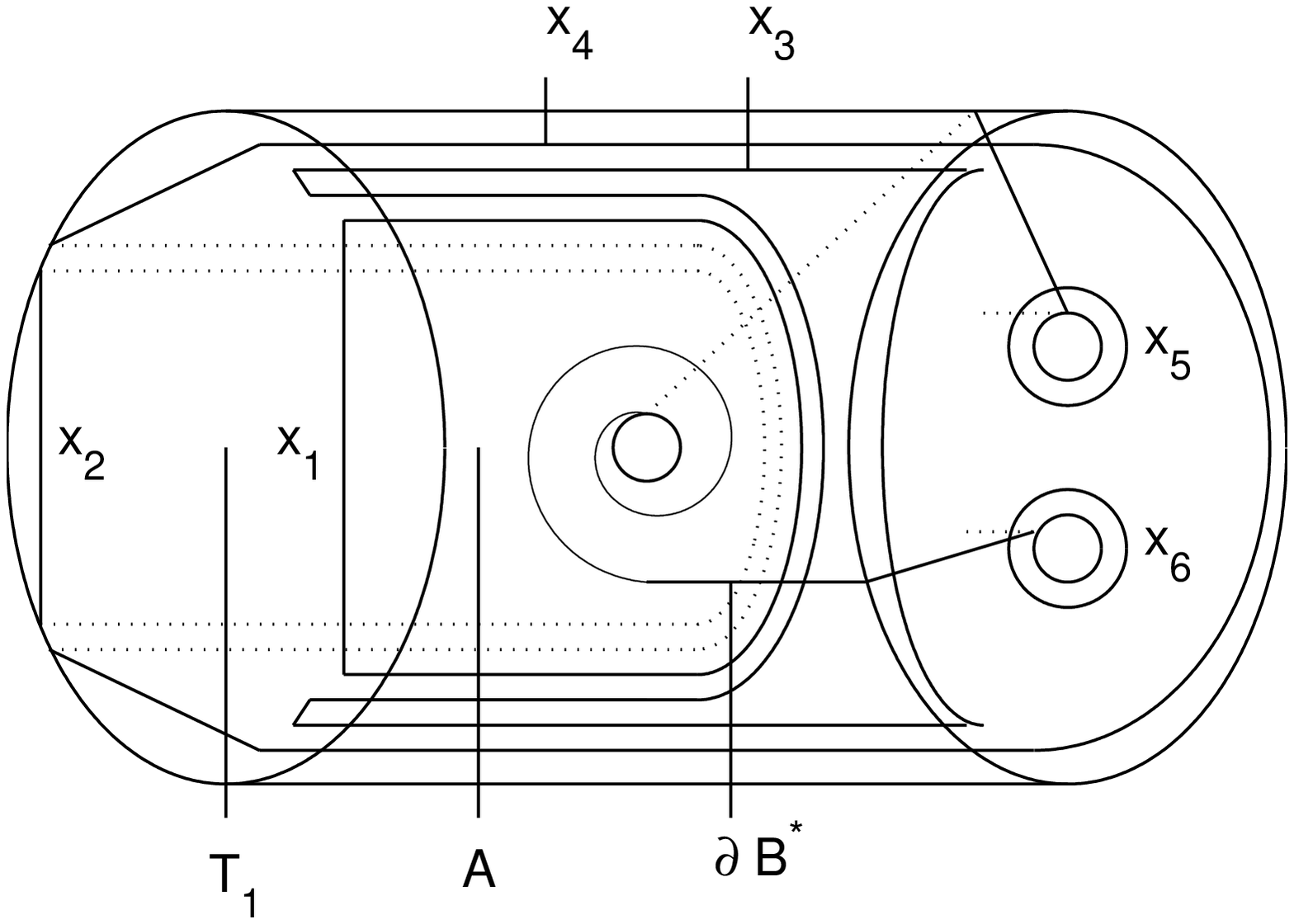}
\begin{center}
Figure 5.14
\end{center}
\end{center}
\begin{center}
\includegraphics[totalheight=5cm]{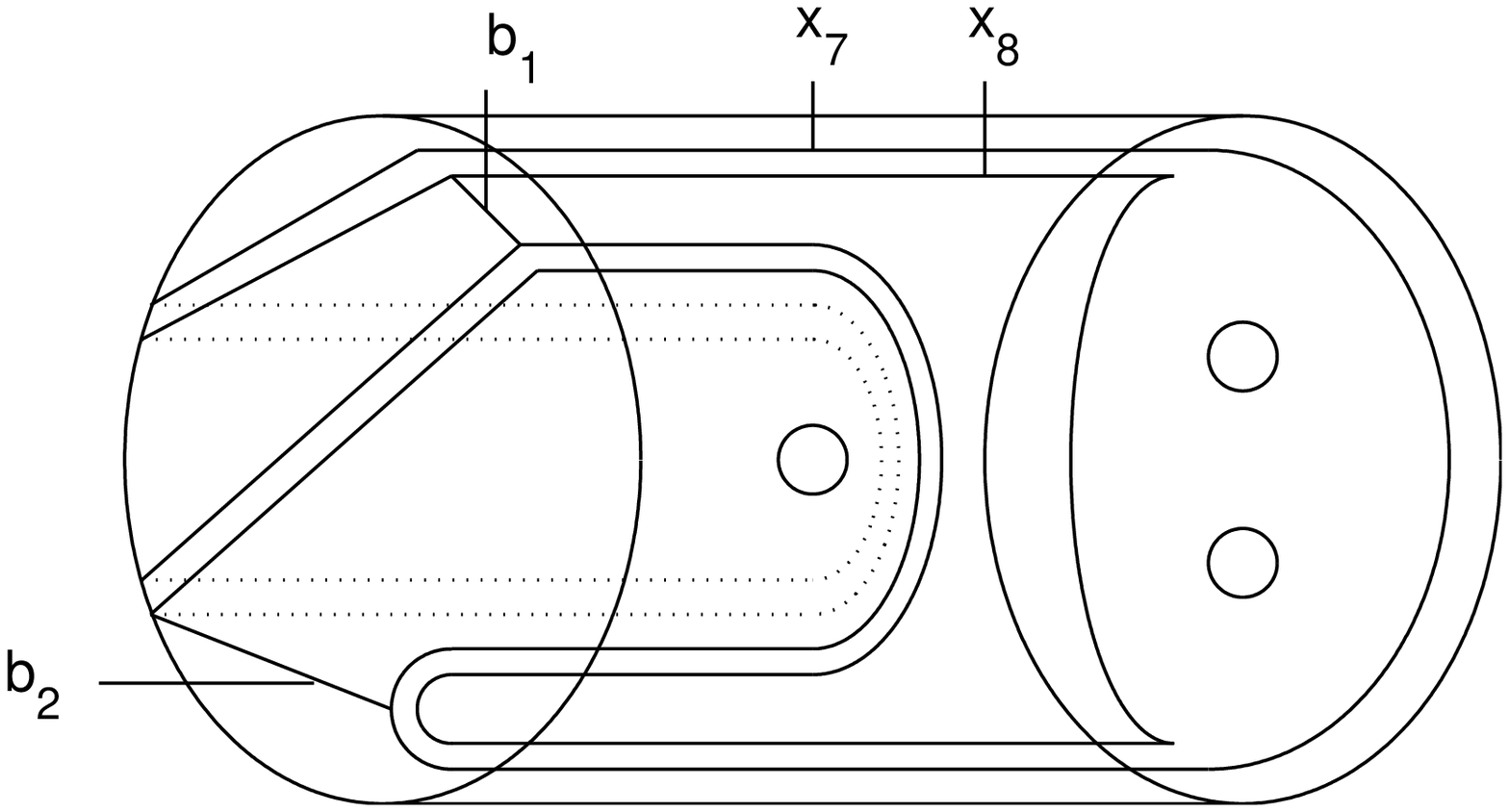}
\begin{center}
Figure 5.15
\end{center}
\end{center}

{\bf Lemma 5.5} $H_{K}$ contains no closed essential surface.

{\bf Proof.} \ Suppose, otherwise, that $H_{K}$ contains an
essential closed surface $F$ such that the complexity $C(F)$ is
minimal among all surfaces isotopic to $F$. By Lemma 5.1, the
pattern of $F\cap F_{2}$ is as in one of Figure 5.3(a) and (b).
Furthermore, $\nu(m_{2})=\nu(m_{3})=0$ for any case. By Lemmas 5.3
and 5.4, the pattern of $F\cap F_{1}$ is as in one of Figures 5.7
and 5.11. Furthermore,  $\nu(n_{2})=\nu(n_{3})=\nu(n_{4})$ for any
case. By Lemma 5.0, $\nu(n_{1})+\nu(n_{2})=\nu(m_{1})$.

In $M_{2}$,  the pattern of $F\cap F_{1}$ can be labeled as in one
of in Figure 5.13(b) or (c), and the pattern of $F\cap F_{2}$ can
be labeled as in Figure 5.13(a).

Note that $W_{2}, W_{4}, W_{8}, W_{10}$ separate $M_{2}$ into four
solid tori $J^{1},J^{2},J^{4},J^{5}$ and a handlebody of genus two
$H^{'}$ such that $A_{2i}\subset J^{i}$ for $i=1,2,4,5$ and
$A_{6}\subset H^{'}$. Let $S=F\cap H^{'}$.

Now we claim that $\nu(n_{2})=\nu(n_{3})=\nu(n_{4})=0$. There are
two cases:

Case 1. \ the pattern of $F\cap F_{1}$ is as in Figure 5.13(b).

Now each component of $\partial S$ is contained in one of the
eight families $x_{1},\ldots,x_{8}$ as in Figures 5.14 and 5.15
where the boundary components of $\partial S$ contained in
$\cup_{i=1}^{4} x_{i}$ are produced by cutting along the arcs in
$F\cap (W_{2}\cup W_{4}\cup W_{8}\cup W_{10})$ whose endpoints lie
in $m_{1}\cup n_{1}$ and the components of $\partial S$ contained
in $x_{7}\cup x_{8}$ are produced by cutting along the arcs whose
endpoints lie in $n_{2}\cup n_{3}\cup n_{4}\cup m_{1}$, and each
component in $x_{5}\cup x_{6}$ is isotopic to one component of
$\partial A_{6}$. Note that each component lying in  $x_{3}\cup
x_{4}$ is trivial in $\partial H^{'}$. By observation, there are
two disks $D^{1}$ and $D^{2}$ in $\partial H^{'}$ such that
$\partial D^{i}=b_{i}\cup b_{i}^{'}$ where $b_{i}\subset F_{1},
b_{i}^{'}\subset S$ as in Figure 5.15. Back to $M_{2}$, $D^{1}$
and $D^{2}$ are as in Figure 5.7. Thus by doing surgeries on $F$
along $D^{1}$ and $D^{2}$, we can obtain a surface $F^{'}$
isotopic to $F$ such that $|F^{'}\cap W|=|F\cap W|$, $|F^{'}\cap
F_{2}|=|F\cap F_{2}|$, $|(F^{'}\cap M_{1}-X(F^{'}))\cap
W^{'}|<|(F\cap M_{1}-X(F))\cap W^{'}|$, a contradiction.

Case 2. \ the pattern of $F\cap F_{1}$ is as in Figure 5.13(c).

This case is similar to Case 1.

\begin{center}
\includegraphics[totalheight=5cm]{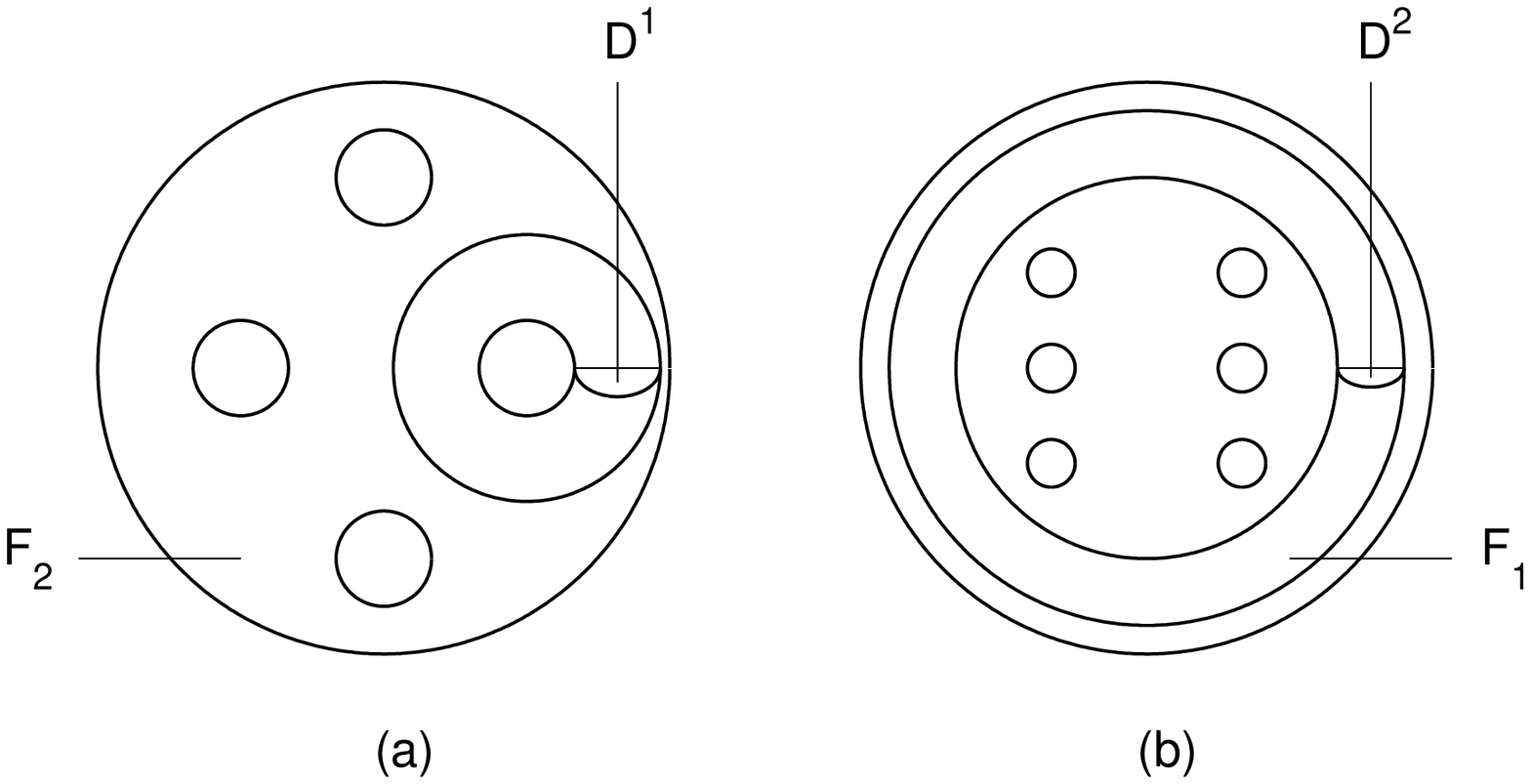}
\begin{center}
Figure 5.16
\end{center}
\end{center}

Now $\nu(n_{2})=\nu(n_{3})=\nu(n_{4})=0$ and $\partial S$ is as in
Figure 5.14.    By construction, there is a disk $B^{*}=H^{'}\cap
D_{6*}$ in $H^{'}$ such that $\partial B^{*}$ intersects each
component in $x_{1}\cup x_{2}\cup x_{5}\cup x_{6}$ in only one
point as in Figure 5.14. Thus $S\cap B^{*}$ offers a
$\partial$-compressing disk $D^{*}$ of $S$ such that $D^{*}$ is
disjoint from $intA_{6}$. We denote by $A$ the annulus bounded by
an outermost component of $x_{1}$, say $e_{1}$, and an outermost
component of $x_{2}$, say $e_{2}$, in $\partial H^{'}$, and
$T_{1}$ the punctured torus bounded by an outermost component of
$x_{1}$ and an outermost component of $x_{2}$ in $\partial
H^{'}$as in Figure 5.14. Now if $\partial D^{*}\cap
\partial H^{'}=a\subset A$, then $e_{1}\cup e_{2}$
bounds an annulus in $S$ parallel to $A$. This means that one
component of $F\cap M_{2}$ is parallel to $\partial H\cap M_{2}$.
Now let $X_{0}(F)$ be a union of components in $F\cap M_{2}$
parallel to $\partial H\cap M_{2}$ or $A_{6}$.

Let $S=(F\cap M_{2}-X_{0}(F))\cap H^{'}$. Then $(F\cap
M_{2}-X_{0}(F))\cap H^{'}\cap B^{*}$ offers a
$\partial$-compressing disk, also denoted by  $D^{*}$, of $S$ such
that $\partial D^{*}\cap
\partial H^{'}=a$.

Now we claim each component of $S$ is isotopic to one component of
$\partial A_{6}$. There are five possibilities:


(1) The two endpoints of $a$ lies in $x_{5}(x_6)$. Then $D^{*}$
can be  moved to be $D^{1}$ as in Figure 5.16(a), Thus by doing a
surgery on $F$ along  $D^{1}$, we can obtain a surface $F^{'}$
isotopic to $F$ such that $|F^{'}\cap W|=|F\cap W|$, $|F^{'}\cap
F_{2}|<|F\cap F_{2}|$, a contradiction.

(2) The two endpoints of $a$ lies in $x_{1}(x_2)$. Then $D^{*}$
can be  moved to be $D^{2}$ as in Figure 5.16(b), contradicting
the minimality of $|F\cap W|$.

(3) One endpoint of $a$ lies in $x_{5}$ and the other lies in $
x_{6}$. Since $\partial B^{*}$ intersects $\cup_{i=1}^{6} x_{i}$
in the order of $x_{6}, x_{3}, x_{1},  x_{2}, x_{4}, x_{5}$, by
the argument in (1), there is an outermost component of $S\cap
B^{*}$ in $B^{*}$, say $b$, which, with an arc $b^{*}$ in
$\partial H^{'}$, bounds an outermost disk $D$ such that
$\partial_{1} b\subset x_{5}$, $\partial_{2}b \subset x_{6}$ and
$b^{*}$ intersects $A_{6}$ in an arc. Since $S$ is incompressible,
by the standard argument,  the component of $S$ containing $b$ is
parallel to $A_{6}$, a contradiction.

(4) One endpoint of $a$ lies in $x_{1}$ and the other lies in $
x_{2}$. Then $\partial_{1} a\subset c_{1}$ and  $\partial_{2}
a\subset c_{2}$, where $c_{1}$ is a component of $x_{1}$, $c_{2}$
is a component of $x_{2}$. We denote also by $A$ the annulus
bounded by $c_{1}, c_{2}$ in $\partial H^{'}$ and $T_{1}$ the
punctured torus bounded by $c_{1}, c_{2}$ in $\partial H^{'}$.
Note that $a$ is disjoint from $intA_{6}$ and $A_{6}\subset
T_{1}$. Hence $a\subset A$. By the above argument, the component
of $F\cap M_{2}$ consisting of $c_{1}$ and  $c_{2}$  is parallel
to $\partial H\cap M_{2}$. By the definition of $S$,  this is
impossible.

(5) One endpoint of $a$ lies in $x_{1}\cup x_{2}$ and the other
lies in $x_{5}\cup x_{6}$.

Since $S$ is incompressible, each component $c$ of $x_{3}\cup
x_{4}$ bounds a disk $D_{c}$ in $S$ parallel to a disk $D_{c}^{*}$
on $\partial H^{'}$, see Figure 5.14. Let $S^{*}=S-\cup_{c\in
x_{3}\cup x_{4}} D_{c}$. Note that $\partial B^{*}$ intersects
$\cup_{i=1}^{6} x_{i}$ in the order: $x_{6}, x_{3}, x_{1}, x_{2},
x_{4}, x_{5}$. Hence each component of $S\cap B^{*}$ is an arc $b$
such that $\partial_{1} b\subset x_{1}\cup x_{2}$ and
$\partial_{2}b\subset x_{5}\cup x_{6}$. Otherwise, there is an
outermost component $b^{*}$ of $S^{*}\cap B^{*}$ in $B^{*}$ such
that $\partial b^{*}$ is as in one of the above four cases, a
contradiction.

Now each component of $S\cap B^{*}$ is an arc $b$ such that
$\partial_{1} b\subset x_{1}\cup x_{2}$ and $\partial_{2}\subset
x_{5}\cup x_{6}$. Let $H^{*}=H^{'}-B^{*}\times (0,1)$ and
$S^{**}=S^{*}-B^{*}\times (0,1)$, where $B^{*}\times I$ is a
regular neighborhood of $B^{*}$ in $H^{'}$. Then $H^{*}$ is a
solid torus. Since each component of $x_{1}\cup x_{2}\cup
x_{5}\cup x_{6}$ intersects $\partial B^{*}$ in one point, each
component $h$ of $\partial S^{**}$ is obtained by doing a band sum
of one component $h_{1}$ of $x_{5}\cup x_{6}$ and one component
$h_{2}$ of $x_{1}\cup x_{2}$ along a component of $S^{*}\cap
B^{*}$.  Since $h_{1}= 1\in\pi_{1}(H)$, $h_{2}\neq 1\in
\pi_{1}(H)$, so $h\neq 1\in\pi_{1}(H^{*})$. Recall the disk
$B_{2}$ in $H$ defined in Section 2. Now $B_{2}\cap H^{'}$ is a
planar surface $P$ such that one component of $\partial P$, say
$\partial_{1} P$, is disjoint from $A_{6}$ and the other
components of $\partial P$ lie in $intA_{6}$. Furthermore,
$\partial_{1} P$ intersects each component in $x_{1}\cup x_{2}$ in
one point. Hence $P-B^{*}\times (0,1)$ is a properly embedded disk
in $H^{*}$ which intersects each component of $\partial S^{**}$ in
one point. This means that each component of $S^{**}$ is an
annulus $A$ parallel to each component of $\partial H^{*}-\partial
A$.

Suppose that $D$ is a $\partial$-compressing disk of $A$ in
$H^{*}$ such that the arc $\alpha=D\cap\partial H^{*}$ lies on the
annulus $A^{*}$ on $\partial H^{*}$ which contains the disk
$A_{6}-B^{*}\times (0,1)$. Then $D$ is disjoint from $x_{3}\cup
x_{4}$.  Since the disk $D^{*}=B^{*}\times\bigl\{0,1\bigr\}\cup
(A_{6}-B^{*}\times (0,1))$ intersects $\partial A^{*}$ in two
arcs, $D$ can be moved to have the arc $\alpha$ lying on
$A^{*}-D^{*}$. Furthermore, since each component $h$ of $\partial
S^{**}$ is obtained by doing a band sum of one component $h_{1}$
of $x_{5}\cup x_{6}$ and one component $h_{2}$ of $x_{1}\cup
x_{2}$, we may assume that $\partial \alpha\subset x_{1}\cup
x_{2}$. Hence $D$ is also a $\partial$-compressing disk of $S^{*}$
in $H^{'}$. By the above argument, this is impossible.

By the above argument, if one component of $F\cap (F_{1}\cup
F_{2})$ is  parallel to $\partial E_{1}$ or $\partial E_{2}$ then
it is parallel to $\partial H$.  Suppose that each component of
$F\cap (F_{1}\cup F_{2})$  is isotopic to one component of
$\partial A_{i}$. By the minimality of $C(F)$, $F$ is disjoint
from $W_{i}$ for $i\neq 6$ and  $F$ is also disjoint from
$\overline{\partial N(B^{*}\cup A_{6})-\partial H^{'}}$ in
$H^{'}$. Thus each component of $F\cap M_{j}$ is an annulus
parallel to $A_{i}$ for some $i$. That means that $F$ is isotopic
to $T$, a contradiction. \qquad Q.E.D.

{\bf The proof of  Proposition 3.0.} \  The proof follows
immediately from Lemma 4.1, 4.3, 4.4 and 5.5 and Theorem 1 in
[SW].\qquad Q.E.D.

{\bf Acknowledgement:} \ The authors thank Professor Fengchun Lei
for some helpful discussions. The authors would like to express
their thanks to the referee for the careful reading of the paper
and pointing out a mistake of the previous version of the paper.

\vskip 0.5 true cm

{\bf References.}

\vskip 0.4 true cm

[CGLS] \ M. Culler, C. Gordon, J. Luecke, and P. Shalen, Dehn
surgery on knots, Ann. of Math., 125(1987), 237-300.

[Go] \ C. Gordon, Combinatorial methods in Dehn surgery, Lectures
at knots, Tokyo, 1996.

[Ha] \ A. Hatcher, On the boundary curves of incompressible
surfaces. Pacific J. Math. 99 (1982), 373-377.

[Jo] \ W. Jaco, Adding a 2-handle to a 3-manifold, An application
to Property R, Proc. Amer. Math. Soc., 92(1984), 288-292.

[L] \ M. Lackenby, Attaching handlebody to 3-manifolds, Geometry
and
                 Topology, Vol. 6 (2002), 889-904 (2002)

[Q] \ R.F. Qiu,  Incompressible surfaces in handlebodies and
closed 3-manifolds of Heegaard genus two, Proc.Amer. Math. Soc.,
128(2000), 3091-3097.

[QW1] \ R.F. Qiu and S.C. Wang,  Small knots and large handle
additions, Comm. Anal. Geom. 13(2005).

[QW2] \ R.F. Qiu and S.C. Wang, Simple, small knots in
handlebodies, Topology and Its applications, 144(2004), 211-227.

[SW] \ M. Scharlemann and Y. Wu, Hyperbolic manifolds and
degenerating handle additions. J. Aust. Math. Soc. (Series A) 55
     (1993), 72-89.

[T] \ W. Thurston, Three dimensional manifolds, Kleinian groups
and hyperbolic geometry. Bull. AMS, Vol. 6, (1982) 357-388.

\vskip 0.5 true cm

Ruifeng Qiu, Department of Mathematics, Dalian University of
Technology, Dalian, 116023, China

\vskip 0.4cm

Shicheng Wang, Institute of Mathematics, Peking University
Beijing, 100871, China

\end{document}